\documentclass[11pt]{amsart}

\usepackage[T1,T5]{fontenc}
\usepackage[utf8]{inputenc}

\usepackage[left=2cm,right=2cm,top=3cm,bottom=3cm]{geometry}

\usepackage{amsthm,amscd,amssymb}
\usepackage{amsmath}
\usepackage{bbm}
\usepackage{mathrsfs}
\usepackage{rotating}
\usepackage{hyperref}
\usepackage{todonotes}
\usepackage{xcolor}

\def\f2{{\mathbb F}_{2}}
\def\z2{{\mathbb Z}/2}

\def\Q{{\mathcal{Q}}}

\def\D{{\mathcal{D}}}
\def\R{{\mathcal{R}}}

\DeclareMathOperator{\gr}{\mathfrak{gr}}

\DeclareMathOperator{\tr}{Tr}

\DeclareMathOperator{\Span}{Span}
\DeclareMathOperator{\F}{\mathbb{F}}

\DeclareMathOperator{\Map}{Map}

\newcommand{\bmt}[1]{\begin{bmatrix}#1\end{bmatrix}}

\newtheorem{theorem}{Theorem}[section]
\newtheorem{lemma}[theorem]{Lemma}
\newtheorem{proposition}[theorem]{Proposition}
\newtheorem{corollary}[theorem]{Corollary}
\newtheorem{conjecture}[theorem]{Conjecture}

\theoremstyle{definition}
\newtheorem{definition}[theorem]{Definition}
\newtheorem{example}[theorem]{Example}
\newtheorem{remark}[theorem]{Remark}

\newtheorem{warning}[theorem]{Warning}
\numberwithin{equation}{section}

\title[Invariants of truncated ring]{On Modular Invariants of the Truncated Polynomial Rings in low ranks}

  \author[Ha]{Le Minh Ha}
  \address{Faculty of Mathematics-Mechanics-Informatics, University of Science, Vietnam National University, Hanoi}
  \author[Hai]{Nguyen Dang Ho Hai}
  \address{Department of Mathematics, College of Sciences, University of Hue, Vietnam}
 \author[Nghia]{Nguyen Van Nghia}
  \address{Faculty of Natural Sciences, Hung Vuong University, Phu Tho, Vietnam}
  \email{leminhha@hus.edu.vn,ndhhai@husc.edu.vn,nguyenvannghia@hvu.edu.vn}

\date{\today}

\begin{document}
	
	\begin{abstract}
		We verify the conjectures due to  Lewis, Reiner, and Stanton about the Hilbert series of the invariant ring of the truncated polynomial ring for all parabolic subgroups up to rank $3$. This is done by constructing an explicit set of generators for each invariant ring in question. We also propose a conjecture concerning the action of the Steenrod algebra and the Dickson algebra on a certain naturally occurring filtration of the invariant ring under the action of the general linear group.  
	\end{abstract} 
	\maketitle
	\section{Introduction}
	Let $\mathbb{F}_q$ denote the finite field of $q$ elements, where $q$ is a power of a fixed prime number $p$. Let $I_m(n) = (x_1^{q^m}, \ldots, x_n^{q^m})$ denote the $m$-th Frobenius power ideal in the polynomial algebra $S(n) = \mathbb{F}_q [x_1, \ldots, x_n]$. Denote by $\Q_m(n)$ the truncated polynomial ring $S(n)/I_m(n)$. The general linear group $G_n=GL_n(\mathbb{F}_q)$ acts on $S(n)$ by linear substitutions of variables. Since $I_m(n)$ is a $G_n$-invariant ideal, $\Q_m(n)$ inherits a $G_n$-action from $S(n)$. As we will be working with a fixed number of generators $n$, we will usually omit $n$ from the notation.   
    
    This paper addresses the problem of computing invariants of the truncated ring $\Q_m (n)$ under the action of various parabolic subgroups of $G_n$. This problem arises in several contexts, including algebraic combinatorics (\cite{ReinerStantonWhite2004}, \cite{Reiner-Stanton2010}, \cite{LewisReinerStanton2017}, \cite{Drescher-Shepler-2020}), algebraic topology (\cite{Boardman93}, \cite{Hewett96}, \cite{Kuhn1987}) and number theory \cite{Deng23}. In particular, as shown in \cite[Conjecture 1.3]{LewisReinerStanton2017}, knowing the Hilbert series of $\Q_m(n)^{G}$ provides the Hilbert series of the cofixed space $S_G$ of the polynomial ring under the action of the general linear group. Surprisingly, while the invariant ring $S^G$ has been completely determined since the beginning of the 20th century, by the fundamental work of Dickson \cite{Dickson11}, the problem of determining $S_G$ is still an open problem in modular invariant theory.  
 
For a positive integer $n$, a composition of $n$ of length $\ell$ is an ordered set $\alpha = (\alpha_1, \ldots, \alpha_{\ell})$ of positive integers $\alpha_i$ such that $|\alpha| = \alpha_1 + \ldots + \alpha_{\ell}=n$. We will also consider weak compositions, which are ordered sets of non-negative integers. The set of weak compositions of a given length is partially ordered by declaring that $\beta \leq \alpha$ iff $\beta_i \leq \alpha_i$ for all $i$. If $\alpha$ is a composition of $n$, we denote by $P(\alpha)$ the corresponding parabolic subgroup of $G_n$. In \cite{LewisReinerStanton2017}, Lewis, Reiner, and Stanton proposed a marvelous set of conjectures about the Hilbert series of the invariant ring $\Q_m(n)^{P(\alpha)}$ for arbitrary $m$ and any composition $\alpha$ of $n$, expressed in terms of a new combinatorial object called the $(q,t)$-binomial coefficient, introduced in \cite{Reiner-Stanton2010}:  
 
	\begin{conjecture}(\cite[Parabolic Conjecture 1.5]{LewisReinerStanton2017})\label{conj: parabolic} The Hilbert series of the ring of invariants of $\Q_m(n)$ under the action of the parabolic subgroup $P(\alpha)$ is: 
		\[
	C_{\alpha,m} (t)= \sum_{\beta, \beta \leq \alpha, |\beta| \leq m} t^{e(m,\alpha,\beta)}  \;  \begin{bmatrix}
			m \\\beta, m- |\beta|
		\end{bmatrix}_{q,t}, 
		\]
	where $e(m, \alpha, \beta) = \sum (\alpha_i- \beta_i) (q^m - q^{B_i})$, $B_i = \beta_1 + \ldots + \beta_i$, and for a composition $\alpha = (\alpha_1, \ldots, \alpha_{\ell})$ of $d$ with partial sum $\alpha_1 + \ldots + \alpha_i = A_i$, the $(q,t)$ binomial coefficient $\begin{bmatrix}
		d \\ \alpha 
	\end{bmatrix}_{q,t}$ is defined as the quotient 
	\[
	\begin{bmatrix}
		d \\ \alpha 
	\end{bmatrix}_{q,t}  =  \frac{\prod_{j=0}^{d-1} (1-t^{q^d-q^j})}{\prod_{i=1}^{\ell} \prod_{j=0}^{\alpha_i} (1- t^{q^{A_i} -q^{A_{i-1}+j}})}.
	\] 
 \end{conjecture}
 In particular, in the most interesting case where $\alpha = (n)$, they conjectured that 
 \begin{conjecture}(\cite[Conjecture 1.2]{LewisReinerStanton2017}) The Hilbert series of the ring of invariants of $\Q_m(n)$ under the action of the general linear group $G_n$ is: 
		\[
	C_{n,m} (t) = \sum_{k=0}^{\min (n,m)} t^{(n-k)(q^m-q^k)}  \;  \begin{bmatrix}
			m \\ k
		\end{bmatrix}_{q,t}.
		\]
 \end{conjecture}
These conjectures are inspired by properties of the $q$-Catalan and $ q$-Fuss-Catalan numbers, connecting Hilbert series of certain invariant subspaces with the representation theory of rational Cherednik algebras for Coxeter and complex reflection groups. Some partial results about these two conjectures have been obtained (see Drescher-Shepler \cite{Drescher-Shepler-2020}, Goyal \cite{Goyal18}, and the original paper \cite{LewisReinerStanton2017}). Nevertheless, the role of these combinatorial objects in the series remains quite mysterious.  

In our previous work \cite{HHN23}, we proved their parabolic conjecture for the minimal parabolic subgroup $B$ by constructing an explicit basis for the $\mathbb{F}_q$-vector space $\Q_m(n)^B$ for all $m$ and $n$. We also refined their conjectures by proposing an explicit $\mathbb{F}_q$-basis for $\Q_m(n)^{P(\alpha)}$. 
\begin{conjecture}\label{refined conj para} A basis for the space of $P(\alpha)$-invariants of $\Q_m(n)$ is given by the set $\mathcal{B}_m(\alpha)$ consisting of elements of the form
\[
\delta_{B_1+1;m}^{\alpha_1-\beta_1}\bigg(f_1\delta_{B_2+1;m}^{\alpha_2-\beta_2}\big(\cdots f_{\ell-1} \delta_{B_\ell+1;m}^{\alpha_\ell-\beta_\ell}(f_\ell )\big)\bigg)
\]
where $\beta \leq \alpha, |\beta| \leq m$, $B_i = \beta_1 + \ldots + \beta_i$ (by convention, $B_0 =0$) and $f_i \in \Phi^{B_{i-1}} \Delta^{m-B_{i-1}}_{\beta_i} \subseteq \Delta^{m-B_{i}}_{B_i}$. \end{conjecture}  

Here $\delta_{a;b}$ is the operator constructed in \cite{HHN23} and $\Delta^m_s$ is a certain subspace of the Dickson algebra $\mathcal{D}_s$ (\cite{Dickson11}) of all $GL_s$-invariant polynomials:
\[
\Delta^m_s \subset \mathcal{D}_s = \mathbb{F}_q [Q_{s,s-1}, \ldots, Q_{s,1}, Q_{s,0}]. 
\]
The Dickson invariants $Q_{s,i}$ are crucial to our investigation and will be reviewed in detail in Section~\ref{sec: Recollection}. The Frobenius-like operator $\Phi$ is a ring map from $\mathcal D_s$ to $\mathcal D_{s+1}$ which sends $Q_{s,i}$ to $Q_{s+1,i+1}$. 

The reader of \cite{HHN23} may notice that in this conjecture, we describe a member of the conjectural generating set $\mathcal{B}_m(\alpha)$ explicitly rather than by induction. In addition, we make use of the core subspace $\Delta^m_s$ (see \ref{subsec:Delta subspace} below) of the Dickson algebra instead of $\nabla^m_s$ in the original conjecture. The two candidates for \emph{core spaces} can be used interchangeably, but $\Delta^m_s$ is more convenient for explicit calculation, as we will see later. 


In the case $\alpha = (n)$ so that $P(\alpha)$ is the full general linear group, we have 
\begin{conjecture}\label{refined conj G}
	 The set $\mathcal{B}_{m}(n)$ consisting of the following family of elements:  
	 \[
	 \delta_{s+1;m}^{n-s} (f), f \in \Delta^m_s, 0 \leq s \leq \min (m,n),
	 \]
	 forms a basis for the $\mathbb{F}_q$-vector space of $G_n$-invariants of $\Q_{m}(n)$.  
\end{conjecture} 
 It is a simple counting problem that our proposed bases have the correct Hilbert series as predicted. Hence our conjectures refine the two conjectures proposed by Lewis, Reiner, and Stanton. Furthermore, they also explain where the summands of the Hilbert series occur in the original conjectures: The summands of $C_{n,m} (t)$ are just the Hilbert series of $\Delta^m_s$, $0 \leq s \leq \min (m,n)$, shifted by an appropriate degree due to the iterated delta operator.  

 In \cite{HHN23}, we have verified the refined parabolic conjecture \ref{refined conj para} for the Borel subgroup, that is when $\alpha = 1^n$. The aim of this paper is to verify our conjectures for other parabolic subgroups, including the full general linear group, for $n \leq 3$. As an immediate corollary, the parabolic conjecture of Lewis, Reiner and Stanton is true in these cases. 
	
	\begin{theorem}\label{thm: main} 
		The refined parabolic conjectures \ref{refined conj para} and \ref{refined conj G}  are true for $n \leq 3$ and all $m \geq 1$. Hence, the parabolic conjecture of Lewis, Reiner, and Stanton is true in these cases.  
	\end{theorem}
	
 In this paper, we choose to work as explicitly as possible to illuminate the underlying mechanism of our construction and to highlight the additional structure present in our basis. Our computations reveal an interesting property of the invariant rings when considered as modules over the Dickson algebra. We also describe the action of the Steenrod algebra. In future work, we hope to revisit this problem in greater generality. 
 
  When $n=2$, the only parabolic subgroup other than the Borel subgroup is $G_2$ itself, and the parabolic conjecture has already been verified in the original paper \cite{LewisReinerStanton2017}. However, the authors worked with the dual situation involving cofixed spaces and used a specific property of the bivariate case to avoid the truncated ring. As a result, it remains unclear whether the invariant subspace of $\Q_{m}(n)$, even in rank $n=2$, can be described explicitly using their method.     

  In \cite[Corollary 4.3]{Goyal18}, Goyal constructed a family of $G_n$-invariant polynomials in $\Q_m(n)$ which corresponds to the family $\delta_2^{n-1} (\Delta^m_1)$ in our notation. 
  

  To simplify the notation, we will often drop the index $n$ in $\Q_m (n), S(n)$, and $G_n$ when the number of variables is well-understood from the context. If the Frobenius power $m$ is fixed, we write $\delta_s$ instead of $\delta_{s;m}$. It is also convenient to write  $(e_1, \ldots, e_s) \in \Delta^m_s$ to indicate a Dickson monomial $Q_{s,s-1}^{e_1} \ldots Q_{s,0}^{e_s}$ in $\Delta^m_s$. For any nonnegative integer $a$, we denote by $[a]_q$ the $q$-integer $(q^a-1)/(q-1)$ with the convention that $[0]_q=1$.   
    
  Here is our first result in rank $2$.  
 		\begin{proposition}\label{prop:GL2}  
 			The set $\mathcal{B}_{m}(2)$ consisting of 3 families below forms a basis for the invariant subspace $\Q_m(2)^G$ of the full general linear group. 
 			\begin{enumerate}
 				\item $\delta_1(\delta_1(1))=x_1^{q^m-1}x_2^{q^m-1}$.  
 				\item $\delta_2 (\Delta^m_1)$.
 				\item $\Delta^m_2$. 
 			\end{enumerate}
    \end{proposition}
  
 In rank 3, other than the Borel subgroups, there are two proper parabolic subgroups. Their corresponding rings of invariants are described in the next two propositions.  
	\begin{proposition}\label{prop:G21} The set $\mathcal{B}_{m}(2,1)$ consisting of 6 families below forms a basis for the invariant subspace $\Q_m^{P(2,1)}$ of the parabolic subgroup corresponding to the composition $(2,1)$.
	 
		\begin{enumerate}
					\item $Q_{2,1}^{i_1} Q_{2,0}^{i_2}Q_{3,2}^i $, $(i_1,i_2) \in \Delta^{m}_2$, $i <[m-2]_q$.    
					\item $Q_{2,1}^{i_1} Q_{2,0}^{i_2}\delta_3 (1) $, $(i_1,i_2) \in \Delta^{m}_2$. \item $\delta_2 (Q_{1,0}^{i_1} Q_{2,1}^{i})$, $i_1 <[m]_q$, $i < [m-1]_q$.   
					\item $\delta_2 (Q_{1,0}^{i_1} \delta_2(1))$, $i_1 < [m]_q$. 
					\item $\delta_1(\delta_1(Q_{1,0}^i))$, $i < [m]_q$. 
					\item $\delta_1(\delta_1(\delta_1(1)))$. 
				\end{enumerate}
	\end{proposition}
	
	\begin{proposition}\label{prop:G12} The set $\mathcal{B}_{m}(1,2)$ consisting of 6 families below forms a basis for the invariant subspace $\Q_m^{P(1,2)}$ of the parabolic subgroup corresponding to the composition $(1,2)$. 
	 
    \begin{enumerate}
					\item $Q_{1,0}^{j_1}  Q_{3,2}^{i_1} Q_{3,1}^{i_2}$, $j_1 < [m]_q, (i_1,i_2) \in \Delta^{m-1}_2$. 
					\item $Q_{1,0}^{j_1} \delta_3 (Q_{2,1}^{j_2})$, $j_1 < [m]_q, j_2 <[m-1]_q$. 
					\item $Q_{1,0}^{j_1} \delta_2(\delta_2(1))$, $j_1 <[m]_q$. 
					\item $\delta_1(Q_{2,1}^{i_1}Q_{2,0}^{i_2})$, $(i_1,i_2) \in \Delta^m_2$. 
					\item $\delta_1(\delta_2(Q_{1,0}^{j_1}))$, $j_1 < [m]_q$. 
					\item $\delta_1(\delta_1(\delta_1(1)))$.    \end{enumerate}
	\end{proposition}

	 Our final result, which is also the most technical, is for the full general linear group in rank 3: 
	
		\begin{proposition}\label{prop:GL3} 		 
			The  set $\mathcal{B}_{m}(3)$ consisting of 4 families below forms a basis for the invariant subspace $\Q_m(3)^{G}$ of the full general linear group. 
		\begin{enumerate}
			\item $\delta_1(\delta_1(\delta_1(1)))=x_1^{q^m-1} x_2^{q^m-1} x_3^{q^m-1}$.  
			\item $\delta_2(\delta_2(\Delta^m_1))$. 
			\item $\delta_3 (\Delta^m_2)$. 
			\item $\Delta^m_3$. 
		\end{enumerate}
	\end{proposition}
	
Since the case of the Borel subgroup, which corresponds to the composition $(1,1,1)$ of $n=3$, was treated in \cite{HHN23}, the propositions stated above together imply the main Theorem \ref{thm: main}. 

It is worth noting that when $m$ is small relative to $n$, all results still hold, provided they are interpreted appropriately. For example, when $n=3$ and $m=2$, the set $\Delta^2_3$ is empty, and $\delta_{3;2} (\Delta^2_2)$ consists of just a single element, namely $\delta_{3;2}(1) = 1 \in \mathbb{F}_q$ in degree $0$. To avoid excessive attention to these special cases (whose proofs are straightforward), we will generally assume that $m$ is sufficiently large compared to $n$ (in fact, $m\ge n$ is enough!).    

In the process of establishing these results, we have also computed the structure of these invariant spaces considered as modules over the Dickson algebra and the mod $q$ Steenrod algebra.

\begin{definition}\label{def:filtration}
   For each $0 \leq k \leq \min (m,n)$, let $\mathcal{F}_{n,k}$ denote the subspace of $\Q_m(n)^G$ spanned by the subsets $\delta_{s+1}^{n-s} (\Delta^m_s)$ where $0 \leq s \leq k$:
\[
\mathcal{F}_{n,k} = \Span \{\delta_{s+1}^{n-s} (f) \colon f \in \Delta^m_s, 0 \leq s \leq k\}.  
\]
\end{definition}
These subspaces form an increasing filtration of $\Q_m(n)^G$, starting with the one-dimensional subspace $\mathcal{F}_{n,0}$, which is spanned by the top-degree class $(x_1 \ldots x_n)^{q^m-1}$. At least in the cases we have computed, this filtration is exhaustive and exhibits remarkable properties. We propose the following:  
\begin{conjecture}\label{conj: filtration} 
For each $1 \leq k < \min (m,n)$,  $\mathcal{F}_{n,k}$ is an $\mathcal{A}$-submodule as well as a $\D_n$-submodule of $\Q_m (n)^G$. Moreover, $\mathcal{F}_{n,k}$ is annihilated by the Dickson invariants $Q_{n,0}, Q_{n,1} \ldots Q_{n,n-k-1}$. 
 \end{conjecture}
This conjecture suggests that the filtration $\mathcal{F}_{n,*}$ may admit an alternative description of a topological nature. We are able to verify this conjecture up to rank $3$: 
\begin{theorem}\label{thm:filtration}
The conjecture \ref{conj: filtration} is true for $n \leq 3$. 
\end{theorem}
   
The layout of this paper is as follows: Section \ref{sec: Recollection} recalls basic results in modular invariant theory and the key constructions from our previous work. In Section \ref{sec: lower bound}, we establish a lower bound on the total dimension of the invariant subspaces. The subsequent four sections are devoted to proving Propositions \ref{prop:GL2} through \ref{prop:GL3}. The final section discusses the actions of the Dickson algebra and the Steenrod algebra, and presents the proof of Theorem \ref{thm:filtration}.

\subsection*{Acknowledgement}  This work was partially completed while the second and third authors were visiting the Vietnam Institute for Advanced Study in Mathematics (VIASM). The third author gratefully acknowledges support from Hung Vuong University’s Fundamental Research Program (Project No. 03/2024). We would like to thank Vo Thanh Tung for providing some computer calculations during the initial stages of this project. We are also grateful to the anonymous referee for their careful reading and numerous suggestions, which significantly improved the exposition. 
	
	\section{Recollection} \label{sec: Recollection}
	In this section, we recall some background materials about modular invariant theory and the main construction in \cite{HHN23}. 
	
	\subsection{Upper triangular invariants and Dickson invariants}
The canonical projection $S \to \Q_m$ is clearly a $G$-equivariant map. Thus, for each parabolic subgroup of $G$,  there is an induced $\mathbb{F}_q$-algebra homomorphism 
	\[
	S^G \subset S^{P} \to \Q_m^{P},
	\]	
 which makes the invariant ring $\Q_m^{P}$ into a module over $S^G$. The structure of the invariant ring $S^G$, which is the Dickson algebra \cite{Dickson11}, is well-known and plays a significant role in our investigation.    
 
 For each positive integer $k$, denote by $V_k$ the product  
	\[
	V_k (x_1, \ldots, x_k) = \prod_{\lambda_i \in \mathbb{F}_q} (x_k + \lambda_1 x_1 + \ldots + \lambda_{k-1} x_{k-1}).  
	\] 
It is well-known that the invariant ring $S^B$ under the Borel subgroup of upper triangular matrices is a polynomial algebra generated by $V_1^{q-1}, V_2^{q-1}, \ldots, V_n^{q-1}$ (see for example, \cite[Theorem 3.4]{Mui75}).  
	Let $X$ be an indeterminate, and recall the fundamental equation
	\[
	V_{n+1} (x_1, \ldots, x_n, X) = X^{q^{n}} + \sum_{i=0}^{n-1} (-1)^{n-i} Q_{n,i} X^{q^i}.
	\]
	The polynomials $Q_{n,i} = Q_{n,i} (x_1, \ldots, x_n)$ are evidently $G$-invariant, and are called Dickson invariants because of the following celebrated theorem due to Dickson: 
	\[
	\D_n = S^G = \mathbb{F}_q [Q_{n,0}, \ldots, Q_{n,n-1}].
	\]
	The Dickson invariants $Q_{n,i}$ can be described more explicitly as follows.  For any $n$-tuple of non-negative integers $(r_1,\dots,r_n)$, put  
	\[
    [r_1,\dots,r_n]=\det(x_i^{q^{r_j}})_{1\leq i,j\leq n}.
    \] 
	In particular, let $L_n = L_n (x_1, \ldots, x_n)$ denote $[0, 1\ldots,n-1]$. We then have  
	\[
	 Q_{n,i}=[0,\ldots,\hat i,\ldots, n]/L_{n}.
	\]
	In particular, $Q_{n,n}=1$ and $Q_{n,0}=L_n^{q-1}$. By convention, we set $Q_{n,i}=0$ if $i<0$ or $i>n$. The following relations are well-known (see for example, \cite[Proposition 1.3]{Wilkerson83}): 
 \begin{proposition}\label{prop: Dickson relation}
 Let $n$ be a positive integer and $0 \leq i \leq n$. The following relations hold: 
 \end{proposition}
 \begin{itemize}
     \item[(i)] $L_n=V_1 V_2 \cdots V_n$. 
     \item[(ii)] $Q_{n,i}= V_n^{q-1}Q_{n-1,i} + Q_{n-1,i-1}^q$.   
 \end{itemize}
Note that $L_n$ is, up to sign, the product of all lines in the $\mathbb{F}_q$-vector space spanned by $ x_1,\ldots,x_n$. 
\subsection{The $\Delta$ subspace of the Dickson algebra}\label{subsec:Delta subspace}	
Lewis, Reiner and Stanton \cite[Subsection 7.4]{LewisReinerStanton2017} speculated that a portion of the Dickson algebra studied in \cite[Section 5]{Reiner-Stanton2010} might form part of an $\mathbb{F}_q$-basis for $\Q_m(n)^G$. Our work confirms this intuition. In fact, this portion, which we call the Delta subspace, plays a central role. We now review its construction and some elementary properties.    

Given a positive integer $s$ and a partition $\lambda_1 \geq \lambda_2 \geq \ldots \geq \lambda_s$ with at most $s$ nonzero parts, we say that a Dickson monomial $Q_{s,s-1}^{e_1} \ldots Q_{s,0}^{e_s}$ in $\D_s$ is of type $(\lambda_1, \ldots, \lambda_s)$ (or q-compatible with $(\lambda_1, \ldots, \lambda_s)$, in the language of \cite[Definition 5.4]{Reiner-Stanton2010}) if for each $1 \leq i \leq s$, 
\[
e_i \in \bigg[\frac{q^{\lambda_i}- q^{\lambda_{i+1}}}{q-1}, \frac{q^{\lambda_i+1}- q^{\lambda_{i+1}}}{q-1}\bigg). 
\]
By convention, $\lambda_{s+1}=0$. Every Dickson monomial in $\D_s$ is $q$-compatible with a unique partition $(\lambda_1, \ldots, \lambda_s)$ of length $\leq s$. 

Now suppose $s \leq m$, and we will restrict our attention to partitions whose Ferrers diagrams fit inside a $s \times (m-s)$ rectangle, that is, partitions with at most $s$ nonzero parts and $m-s \geq \lambda_1$. Define $\Delta_{(\lambda_1, \ldots, \lambda_s)}$ as the subset of the Dickson algebra $\D_s$ consisting of Dickson monomials of type $(\lambda_1, \ldots, \lambda_s)$, and let $\Delta^m_s$ denote the disjoint union of all subsets over all such partitions fitting inside an $s \times (m-s)$ rectangle. By convention, we put $\Delta^m_s = \emptyset$ if $s> m$, and $\Delta^m_0 = \{1\}$. 
  
  The Dickson monomials in $\Delta^m_s$ for $s \leq \min (m,n)$ together with the delta operators recalled below are the building blocks for our basis. For this reason, we will say that Dickson monomials in $\Delta^m_s$ \emph{essential} monomials. We will also work with the \emph{edge} monomials, which are the Dickson monomials not in $\Delta^m_s$ that have the form $Q_{s,i} f$ where $f \in \Delta^m_s$, for some $0 \leq i \leq s-1$. Clearly, the edge monomials are those  such that $e_j = \frac{q^{\lambda_j+1}- q^{\lambda_{j+1}}}{q-1}$ for some $1 \leq j \leq s$ and 
  $e_i \in [\frac{q^{\lambda_i}- q^{\lambda_{i+1}}}{q-1}, \frac{q^{\lambda_i+1}- q^{\lambda_{i+1}}}{q-1})$ for all $i \neq j$.  


	\subsection{The delta operator}\label{subsec:delta}
	The operator $\delta_{a;b}$ introduced in \cite{HHN23} is fundamental to our construction. We now recall its definition and main properties. 
	
	\begin{definition}Let $a,b,c$ be positive integers such that $1\le a\le c+1$. If $f$ is a rational function in $c$ variables, then $\delta_{a;b}(f)$ is a function in $(c+1)$ variables, defined as the quotient 
		$$\delta_{a;b}(f)=\frac{\begin{vmatrix}
				x_1 & \cdots &x_a \\
				x_1^q & \cdots &x_a^q \\
				\vdots & \ddots & \vdots\\
				x_1^{q^{a-2}} & \cdots & x_a^{q^{a-2}}\\
				x_1^{q^b} f(\widehat{x_1},x_2,\ldots,x_{c+1})& \cdots & x_a^{q^b}f(x_1,\ldots,\widehat{x_a},\ldots,x_{c+1})
		\end{vmatrix}}{
			\begin{vmatrix}
				x_1 & \cdots &x_a \\
				x_1^q & \cdots &x_a^q \\
				\vdots & \ddots & \vdots\\
				x_1^{q^{a-1}} & \cdots & x_a^{q^{a-1}}
		\end{vmatrix}}.$$
	\end{definition}
	
	Thus if $b \geq a-1$, $\delta_{a;b}$ increases degree by $q^b-q^{a-1}$. By definition, if we separate the first $(a-1)$ variables and write $f$ in the form $f=\sum g(x_1,\ldots,x_{a-1}) h(x_a,\ldots, x_c)$, then 
    \[ 
    \delta_{a;b}(f)=\sum \delta_{a;b}(g) h(x_{a+1},\ldots,x_{c+1}). 
    \]
	
 
It was shown in greater generality in \cite{HHN23} that $\delta_s(f)$ is generally not a polynomial. Still, in the cases that we are interested in, it is a genuine polynomial. Moreover, $\delta_s(f)$ is rarely a Dickson polynomial but it is a $G_s$-invariant modulo $I_m$.  
We give a precise statement and provide an elementary proof for $\delta_3$ to make the paper self-contained. 
 \begin{proposition}
     \label{prop:delta3} If $f$ is a $G_2$-invariant polynomial, then $\delta_3 (f)$ is also a polynomial and is $G_3$-invariant in $\Q_{m}(3)$. 
 \end{proposition}
 \begin{proof}
 	We first show that $\delta_3 (f)$ is a polynomial. Recall the denominator $L_3 = V_1 V_2 V_3$ is a product of nonzero linear forms in the variables $x_1, x_2, x_3$, so it suffices to prove that the numerator of $\delta_3 (f)$ becomes zero whenever there exists a nontrivial linear relation among the $x_i$, which by symmetry, can be assumed to have the following form: 
    \begin{equation}\label{eq: linear relations among xi}
        x_3 = a_1 x_1 + a_2 x_2, \quad (a_1, a_2) \in \F_q\times \F_q.  
    \end{equation}
    We claim that this implies a similar linear relation among the columns of the numerator of $\delta_3 (f)$, that is, 
 \begin{equation}\label{eq: linear relation among columns of numerator}
      x_3^{q^m} f(x_1, x_2) = a_1 x_1^{q^m} f(x_2, x_3) + a_2 x_2^{q^m} f(x_1, x_3). 
 \end{equation}
 Indeed, this is clear if $a_1=a_2=0$. If the scalars $a_1, a_2$ are both nonzero, then there is an equality
 	\[
 	f(x_2, a_1 x_1 + a_2 x_2) =  f(x_1, a_1 x_1 + a_2 x_2) = f(x_1, x_2),
 	\] 
    since $f$ is a $G_2$-invariant. So the right-hand side of \eqref{eq: linear relation among columns of numerator} becomes 
 	\[
 	(a_1 x_1^{q^m}   + a_2 x_2^{q^m}) f(x_1, x_2) = (a_1 x_1  + a_2 x_2)^{q^m} f(x_1, x_2) = x_3^{q^m} f(x_1,x_2).  
 	\]
Finally, if one of the $a_i$ is zero, say $a_1 = 0$, Equation \eqref{eq: linear relation among columns of numerator} is reduced to 
  \[
 (a_2 x_2)^{q^m} f(x_1, x_2) = a_2 x_2^{q^m} f(x_1, a_2 x_2), 
  \]
   which is again true because $f$ is $G_2$-invariant. We have finished the proof that $\delta_3 (f)$ is a polynomial. 
  
  The next step is to show that $\delta_3 (f)$ is $G_3$-invariant modulo $I_{m}$. By symmetry, it suffices to prove that $\delta_3 (f)$ does not change under the operation that sends $x_1$ to $x_1 + x_2$ and fixes $x_2, x_3$. The denominator $L_3$ certainly satisfies this condition. By direct inspection, the difference $\delta_3 (f) (x_1+x_2, x_2, x_3) - \delta_3 (f) (x_1, x_2, x_3) $ equals
 	\[
 	x_2^{q^m} \times  \frac{ 
 		\begin{vmatrix}
 			x_1 & x_2 & x_3\\x_1^q & x_2^q & x_3^q\\   f(x_2,x_3)-f(x_1+x_2,x_3)&  f(x_1+x_2,x_3)-f(x_1, x_3) & 0  
 	\end{vmatrix}}{L_3}.
 	\] 
  We already know that this is a polynomial. Moreover, $x_2$ occurs as a simple factor in the denominator $L_3$, and the determinant in the numerator vanishes when $x_2 = 0$. It follows that this polynomial is divisible by $x_2^{q^m}$, thus trivial in $\Q_m(3)$.  
 \end{proof}
 
 \begin{example}\label{y family} 
 The $G_1$-polynomials are of the form $f(x) = x^{s(q-1)} = Q_{1,0}^s$, $s \geq 0$. Then $\delta_2 (x^{s(q-1)})$ is a $G_2$-invariant polynomial in $\Q_{m} (2)$. For $s < [m]_q$, this is Goyal's $y_s$ family in \cite{Goyal18} where it was defined as 
 \[
 y_s = x_1^{q^m-q} x_2^{s(q-1)} + x_1^{q^m-q-(q-1)} x_2^{(s+1)(q-1)} + \ldots + x_1^{s(q-1)} x_2^{q^m-q}. 
 \]  
Our construction is also valid for $s \geq [m]_q$. However, it is trivial in $\Q_m(2)$ except when $s = [m]_q+1$, in which case 
\[
\delta_2 (Q_{1,0}^{[m]_q +1}) = - x_1^{q^m-1} x_2^{q^m-1} = - \delta_1^2 (1).
\]
It is easy to check that $y_s$ is a genuine $G_2$-invariant (that is, a Dickson polynomial) whenever $s = \frac{q^i-1}{q-1}$ for some non-negative integer $i$. 
 \end{example}
 \begin{example}\label{a family} 
 For $s \leq [m]_q$, $\delta_2^2 (Q_{1,0}^s)$ is a $G_3$-invariant of degree $2(q^m-q)+s(q-1)$ in $\Q_m(3)$. This family of elements coincides with Goyal's polynomials $a_{m,3,s}$ in \cite[Corollary 4.3]{Goyal18} where it was written in a different form:
 \[
a_{m,3,s} = \sum_{(i_1,i_2,i_3)}  x_1^{i_1 (q-1)} x_2^{i_2 (q-1)} x_3^{i_3 (q-1)},  
 \]
 where the sum is taken over all triples $(i_1,i_2,i_3)$ of non-negative integers such that $i_s < [m]_q$ and $i_1 + i_2 + i_3=2 [m]_q - 2 +s$. Again, our definition also works for $s > [m]_q$: 
 \[
 \delta_2^2 (Q_{1,0}^{[m]_q + 2}) = \delta_1^3 (1) = x_1^{q^m-1} x_2^{q^m-1} x_3^{q^m-1},
 \]
 and is trivial in all other cases.  
 \end{example}
\begin{warning}\label{warning}
   It is not true in general that if $f = g$ in $\Q_m(n)$ then $\delta_s (f) = \delta_s (g)$ in $\Q_m(n)$, even in the case where all expressions involved are polynomials. This is evident from Example \ref{y family} above for $f = Q_{1,0}^{[m]_q+1}$ and $g= 0$. Thus when working with iterated delta operators, for example, $\delta_2^2(f) = \delta_2 (\delta_2 f)$, we cannot take $\delta_2 f$ modulo $I_m$ before applying $\delta_2$ again.  
\end{warning}
   
 \subsection{The delta operator and the Dickson algebra}\label{subsec: delta and Dickson}
For each composition $\alpha$ of $n$, the invariant ring $\Q_m(n)^{P(\alpha)}$ is a module over the Dickson algebra $\D_n$. The next result describes how the delta operator interacts with the Dickson algebra in low ranks. 
\begin{proposition} \label{prop: delta and Dickson} We have the following identities in $\Q_m$. 
 \begin{enumerate}
     \item $Q_{s,0} \delta_s (f) = 0$ for all $f \in S$.
     \item $Q_{2,1} \delta_2 (f) = \delta_2 (Q_{1,0}^q f)$ for all $f \in \D_1$. 
     \item $Q_{3,i} \delta_{3} (f) = \delta_3 (Q_{2,i-1}^q f)$ for $i=1$ or $2$ and for all $f \in \D_2$. 
     \item $Q_{3,2} \delta_{2}^2 (f) = \delta_2^2 (Q_{1,0}^{q^2} f)$ for all $f \in \D_1$.  
     \item $Q_{3,1} \delta_2^2 (f) = 0$ for all $f \in \D_1$. 
 \end{enumerate}
\end{proposition}
\begin{proof}
We will make use frequently of the relation \ref{prop: Dickson relation} relating upper-triangular invariants $V_k$ and Dickson invariants $Q_{n,i}$. 

The first part is clear since $Q_{s,0} = L_s^{q-1}$. Note that the numerator of $\delta_s (f)$ always vanishes in $\Q_m$ by considering the Laplace expansion along the last row. 

\noindent \textbf{Proof of (2).} Since $f \in \mathcal{D}_1$, let us assume that $f = Q_{1,0}^s$ for some $s \geq 0$. Since 
\[
Q_{2,1} = V_2^{q-1}(x_1, x_2) + V_1^{q(q-1)} (x_1) = V_2^{q-1}(x_2, x_1) + V_1^{q(q-1)} (x_2),
\]
the numerator of the difference $Q_{2,1} \delta_2 (f) - \delta_2 (Q_{1,0}^q f)$ has the form 
\[
\begin{vmatrix}
    x_1 & x_2 \\
    x_1^{q^m} V_2 (x_2,x_1)^{q-1} x_2^{s(q-1)} & x_2^{q^m} V_2 (x_1,x_2)^{q-1} x_1^{s(q-1)}
\end{vmatrix}, 
\]
Observe that $x_1 V_2 (x_1, x_2) = x_2 V_2 (x_2, x_1) = L_2 (x_1, x_2)$, it follows that 
\[
Q_{2,1} \delta_2 (f) - \delta_2 (Q_{1,0}^q f) = x_2^{q^m} V_2 (x_1, x_2)^{q-2} x_1^{s(q-1)} - x_1^{q^m} V_2 (x_2, x_1)^{q-2} x_2^{s(q-1)} = 0 \quad \text{in $\Q_m(2)$}. 
\]
\noindent \textbf{Proof of (3).}  We will treat the case $i=1$, as the situation for $i=2$ is similar. Consider the difference $Q_{3,1} \delta_3 (f) - \delta_3 (Q_{2,0}^q f)$. Using the identity $Q_{3,1} = V_3^{q-1} Q_{2,1} + Q_{2,0}^q$, the (3,1)-entry in the last row of the determinant in the numerator can be simplified as 
\[
x_1^{q^m} V_3^{q-1}(x_2,x_3,x_1) Q_{2,1}(x_2,x_3) f(x_2,x_3).  
\]
There is a similar description for other entries in this row. Since 
\[
V_3 (x_2,x_3,x_1) L_2 (x_2,x_3) = L_3 (x_2,x_3,x_1) = L_3 (x_1,x_2,x_3),
\] 
the Laplace expansion along the last row of the determinant in the numerator shows that the resulting polynomial belongs to $I_m$, which is zero in $\Q_m(3)$. 

\noindent \textbf{Proof of (4).} From Proposition \ref{prop: Dickson relation}, we have $Q_{3,2} = V_3^{q-1}+ Q_{2,1}^q$. We
write $Q_{3,2}$ in two different ways: 

\begin{equation*}
    Q_{3,2} = V_3^{q-1}(x_2, x_3, x_1)  + Q_{2,1}^{q} (x_2,x_3) = V_3^{q-1}(x_1, x_3, x_2)  + Q_{2,1}^{q} (x_1,x_3). 
\end{equation*}   
It follows that $Q_{3,2} \delta_2 (\delta_2 f) - \delta_2 (Q_{2,1}^q \delta_2 f)$ equals 
\[
\frac{1}{L_2(x_1,x_2)} \begin{vmatrix}
    x_1 & x_2\\
    x_1^{q^m} V_3^{q-1}(x_2,x_3,x_1) (\delta_2 f)(x_2,x_3) &  x_2^{q^m} V_3^{q-1}(x_1,x_3,x_2) (\delta_2 f)(x_1,x_3)
\end{vmatrix}. 
\]
Since $x_2 V_3 (x_2,x_3,x_1)$ and $x_1 V_3 (x_1, x_3, x_2)$ are divisible by $L_2 (x_1,x_2)$, the above expression is a polynomial in the ideal generated by $(x_1^{q^m}, x_2^{q^m})$. Thus 
\[
Q_{3,2} \delta_2 (\delta_2 f) = \delta_2 (Q_{2,1}^q \delta_2 f) \quad \text{in $\Q_m(3)$}. 
\]
At this point, as noted in Warning \ref{warning}, we cannot apply immediately part (2), replacing $Q_{2,1}^q \delta_2 f$ with $\delta_2 (Q_{1,0}^{q^2}f)$ on the right-hand side. Rather, we make use of the above strategy one more time, writing
\[
Q_{2,1} (x_2,x_3) = V_2^{q-1} (x_3,x_2) + V_1^{q^2(q-1)}(x_3), 
\]
and similarly for $Q_{2,1}(x_1,x_3)$. The difference $\delta_2 (Q_{2,1}^q \delta_2 f) - \delta_2^2 (Q_{1,0}^{q^2}f)$ can be written as the sum 
\begin{multline*}
x_3^{q^m} \frac{x_2^{q^m} x_1 f(x_1) V_2^{q^2-q-1}(x_1,x_3) - x_1^{q^m} x_2 f(x_2) V_2^{q^2-q-1}(x_2,x_3)}{L_2 (x_1,x_2)} - \\ x_2^{q^m} f(x_3) \frac{x_1 x_1^{q^m} V_2^{q^2-q-1}(x_3,x_1) - x_2 x_1^{q^m} V_2^{q^2-q-1}(x_3,x_2)}{L_2(x_1,x_2)}. 
\end{multline*}
Now the first summand is a polynomial (since the numerator is divisible by $x_1$ and vanishes whenever there is a linear relation $x_2=a_1x_1$) and is a multiple of $x_3^{q^m}$. Similarly, the second summand belongs to the ideal $(x_2^{q^m})$.  

\noindent \textbf{Proof of (5).} We first apply the same strategy as part (4), writing
\[
Q_{3,1} = V_3^{q-1} (x_1,x_3,x_2) Q_{2,1} (x_1, x_3) + Q_{2,0}^q (x_3,x_2) =  V_3^{q-1} (x_2,x_3,x_1) Q_{2,1} (x_2, x_3) + Q_{2,0}^q (x_3,x_1),
\]
and obtain an equality in $\Q_m(3)$: 
\[
Q_{3,1} \delta_2 (\delta_2 f)) = \delta_2 (Q_{2,0}^q \delta_2 f). 
\]
The right-hand side, by direct inspection, can be written as 
\begin{multline*}
x_3^{q^m} \frac{x_2^{q^m} x_1^2 f(x_1) L_2^{q^2-q-1}(x_1,x_3)- x_1^{q^m} x_2^2 f(x_2) L_2^{q^2-q-1}(x_2,x_3)}{L_2 (x_1,x_2)} + \\ x_2^{q^m} f(x_3) x_3 \frac{x_1 x_1^{q^m}  L_2^{q^2-q-1}(x_1,x_3)-x_2x_1^{q^m} L_2^{q^2-q-1}(x_2,x_3)}{L_2(x_1,x_2)}. 
\end{multline*}
By the same argument as in part (4), the first summand is divisible by $x_3^{q^m}$ and the second by $x_2^{q^m}$. 
\end{proof}

		\subsection{The Borel subgroup}\label{subsec: Borel}
In \cite{HHN23}, we constructed a basis for the invariant ring of $\Q_m(n)$ under the action of the Borel subgroup for arbitrary $m$ and $n$. In the case $n=3$, the $\mathbb{F}_q$-basis $\mathcal{B}_{m}(3)$ for $\Q_m^B$ consists of $4$ families:  
		\begin{enumerate}
 			\item $x_1^{q^m-1} x_2^{q^m-1} x_3^{j_3 (q-1)}$, $j_3 \leq [m]_q$. 
 			\item $x_1^{q^m-1} x_2^{j_2 (q-1)} Q_{2,1} (x_2, x_3)^{j_3}$, $j_2 < [m]_q, j_3 \leq [m-1]_q$. 
 			\item $x_1^{j_1 (q-1)} \delta_2 (Q_{2,1}^{j_2})$, $ j_1<[m]_q$, $j_2 \leq [m-1]_q$.  
 			\item  $x_1^{j_1(q-1)} Q_{2,1}^{j_2} Q_{3,2}^{j_3}$, $j_1 < [m]_q$, $j_2 <[m-1]_q$, $ j_3 \leq [m-2]_q$. 
 		\end{enumerate}

	\section{A lower bound on the total dimension of the invariant subspaces}\label{sec: lower bound}

 The original Parabolic Conjecture predicted that the Hilbert series of the invariant space $\Q_m(n)^{P(\alpha)}$, where $\alpha$ is a composition of $n$, is given by
 the (finite) polynomial $C_{\alpha,m}(t)$, whose explicit form is stated in Conjecture \ref{conj: parabolic}. In particular, the value of $C_{\alpha,m}(t)$ when $t=1$ is expected to equal the total dimension of the finite-dimensional graded $\mathbb{F}_q$-vector space $\Q_m(n)^{P(\alpha)}$. 

 In this section, we will show that: 
 \begin{proposition}\label{prop: lower bound dim} 
     For each $m,n \geq 1$ and any composition $\alpha$ of $n$, the total dimension of the graded vector space  $\Q_m(n)^{P(\alpha)}$  is at least $C_{\alpha,m}(1)$.   
 \end{proposition} 
 This result is already implicit in \cite[Section 6]{LewisReinerStanton2017} and holds for arbitrary $n \geq 1$. This relatively simple observation has an important practical consequence: suppose we can construct a generating set for the invariant subspace $\Q_m(n)^{P(\alpha)}$ such that, by dimension counting, the associated Hilbert series $C'_{m, \alpha}(t)$  does not exceed $C_{\alpha,m} (t)$ (i.e. $f(t) \leq g(t)$ if and only if $g(t)-f(t)$ is a polynomial with non-negative coefficients). Then, since the reverse inequality $C'_{\alpha,m}(1) \geq C_{\alpha,m}(1)$ holds evaluated at $t=1$, we may conclude that these two series must be identical since both are finite degree polynomials with non-negative coefficients. It follows that our generating set is, in fact, a basis for $\Q_m(n)^{P(\alpha)}$. 

 \vspace{.1in}
\noindent \textbf{Proof of  Proposition \ref{prop: lower bound dim}:} 
  Consider the following ungraded quotient of the polynomial algebra: 
\[
\mathcal{R}=\F_q[x_1,\ldots,x_n]/(x_1^{q^m}-x_1,\ldots, x_n^{q^m}-x_n).
\]
 The action of $G=GL_n$ on $\F_q[x_1,\ldots,x_n]$ descends to this quotient. There is a natural filtration:
 \[
 \{0\}=F_{-1}\subset F_0\subset F_1\subset \cdots \subset F_{n(q^m-1)}=\R,
 \]
 where $F_i$ is the image of the set of all polynomials of degree at most $i$ under the canonical projection  $\F_q[x_1,\ldots,x_n]\twoheadrightarrow \R$. Let $\gr_F\R$ denote the associated graded vector space $\bigoplus_{i\ge 0}F_i/F_{i-1}$ and equip $\gr_F \R$ with a ring structure where the multiplication is induced by the product  $F_i F_j\to F_{i+j}$. The canonical $\F_q$-algebra map 
\[
 \F_q[x_1,\ldots,x_n]\to \R \to \gr_F\R, \quad x_i\mapsto \bar x_i\in F_1/F_0,
\]
induces an isomorphism of $G$-equivariant $\F_q$-algebras:
\[
\Q_m (n) \cong \gr_F\R.
\]
On the other hand, if we apply the functor taking $G$-invariants on each short exact sequence of $G$-modules 
\[
0\to F_{i-1}\to F_{i}\to F_i/F_{i-1}\to 0,
\]
and compute the dimension, we obtain an inequality 
\[
\dim (F_i)^{G}\le  \dim (F_{i-1})^{G}+ \dim (F_i/F_{i-1})^{G}
\]
and so 
\[
\dim \R^{G} \leq \dim (\gr_F\R)^{G} = \dim \Q_m(n)^{G}. 
\]
To compute the $\F_q$-dimension of $\R^{G}$, we extend the field $\F_q$ to a field $\F$ which contains $\F_{q^m}$. Then the evaluation map
\[
f(x_1,\ldots,x_n)\mapsto [(v_1,\ldots,v_n)\in \F_{q^m}^{n}\mapsto f(v_1,\ldots,v_n)\in \F],
\]
induces an isomorphism of $\F$-algebras
\[
\F\otimes_{\F_q}\R \xrightarrow{\cong} 
\Map\big(\F_{q^m}^{n}, \F\big).
\]
This isomorphism is also $G$-equivariant where $G$ acts on the right-hand side via the embedding $G = GL_n(\F_q) \subset GL_n(\F_{q^m}).$ It follows that there is an isomorphism
\[
(\F\otimes_{\F_q}\R)^{G} \cong
\Map\big(\F_{q^m}^{n}/G, \F\big).
\]
The $\F_q$-dimension of $(\F\otimes_{\F_q}\R)^{G}$ is thus equal to the cardinality of the orbit set $\F_{q^m}^n/G$. Finally this cardinality is equal to the sum $\sum_{s=0}^{\min(m,n)}\bmt{m\\ s}_{q}$ where $\bmt{m\\s}_{q}$ is the $q$-binomial coefficient which counts the number of $s$-dimensional subspaces of the $\F_q$-vector space $\F_{q^m}$.

In summary, we have shown that 
\[
\dim_{\F_q} \R^{G} = \dim_{F} (\F \otimes_{\F_q} \R^{G}) = \sum_{s=0}^{\min(m,n)}\bmt{m\\ s}_{q} = C_{m,n}(1) \leq \dim \Q_m(n)^{G}. 
\]
The argument above clearly is not affected if one uses a parabolic subgroup $P_{\alpha}$ instead of $G$. The proposition is proved. \qed 

We end this section with the following observation about our proposed basis $\mathcal{B}_m (\alpha)$. 
\begin{lemma}\label{lem:size}
For each composition $\alpha$ of $n$, the Hilbert series for the $\mathbb{F}_q$-space spanned by $\mathcal{B}_m (\alpha)$ is not greater than $C_{\alpha,m}(t)$. 
\end{lemma}
\begin{proof}
This is almost obvious from the definition of $C_{\alpha,m} (t)$ and the relationship between 
the subset $\Delta^m_s \subset \D_s$ and the (multinomial) $(q,t)$-coefficient. Recall that 
\[
C_{\alpha,m} (t) = \sum_{\beta \leq \alpha, |\beta|\leq m} t^{e(m,\alpha, \beta)} \begin{bmatrix}
			m \\\beta, m- |\beta|
		\end{bmatrix}_{q,t},
\]
where $e(m,\alpha, \beta) = \sum (\alpha_i - \beta_i) (q^m-q^{B_i})$. Our conjectural basis explains precisely how the summands in the above expression occur. Indeed, for each composition
$\beta$, with $\beta \leq \alpha$ and $|\beta| \leq m$, $e(m, \alpha, \beta)$ is exactly the total degree raised by the delta operators in a basis element 
\[
\delta_{B_1+1;m}^{\alpha_1-\beta_1}\bigg(f_1\delta_{B_2+1;m}^{\alpha_2-\beta_2}\big(\cdots f_{\ell-1} \delta_{B_\ell+1;m}^{\alpha_\ell-\beta_\ell}(f_\ell )\big)\bigg),
\]
 and the $(q,t)$-multinomial coefficient $\begin{bmatrix}
			m \\\beta, m- |\beta|
		\end{bmatrix}_{q,t}$ corresponds to the number of choices of  $\ell$-tuple $(f_1, \ldots, f_{\ell})$ where $f_i \in \Phi^{B_{i-1}} \Delta_{\beta_i}^{m-B_{i-1}}$. Note that according to \cite[(7.1)]{Reiner-Stanton2010}, we can write
        \[
\begin{bmatrix}
			m \\\beta, m- |\beta|
		\end{bmatrix}_{q,t}
  =  \begin{bmatrix}
			m \\ \beta_1
		\end{bmatrix}_{q,t} \varphi^{\beta_1} \begin{bmatrix}
			m - B_1 \\ \beta_2 
		\end{bmatrix}_{q,t}  \varphi^{\beta_2} \begin{bmatrix}
			m - B_2 \\ \beta_3 
		\end{bmatrix}_{q,t}  \ldots    
\]
where $\varphi$ is the (genuine) Frobenius operator, which reflects how our operator $\Phi$ affects the degree.  
\end{proof}


 The results in this section simplify our tasks considerably. Indeed, in order to show that $\mathcal{B}_m (\alpha)$ is a basis for $\Q_m(n)^{P_{\alpha}}$, it suffices to verify that it is a generating set.

	\section{Proof of Proposition \ref{prop:GL2}}\label{sec:G2}
As remarked in the previous section, it suffices to show that $\mathcal{B}_{m}(2)$ is a generating set for $\Q_{m}(2)^G$. The proof of Proposition \ref{prop:GL2} is done in two steps. Firstly, we use the relative transfer from the known invariant ring $\Q^B$ of the Borel subgroup $B$ to show that a bigger set $\mathcal{B}'$, consisting of 3 families:
\begin{enumerate}
    \item $\delta_1^2 (1) = x_1^{q^m-1} x_2^{q^m-1}$,
    \item $\delta_2 (Q_{1,0}^s) = y_s$, $0 \leq s < [m]_q$,
    \item $\D_2$, 
\end{enumerate}
is a generating set for $\Q_{m}(2)^G$. Then we show that to generate the invariant ring, it is possible to replace the Dickson algebra $\D_2$ in the third family by the smaller subspace $\Delta^m_2$. 
	   
	\begin{lemma}\label{lem: G2 bigger generating set}
		The set $\mathcal{B}'$ consisting of 3 families above spans $\Q_m(2)^G$.  
	\end{lemma}
	\begin{proof}
 In \cite{HHN23}, we described an $\mathbb{F}_q$-basis $\mathcal{B}_{m}(1,1)$ for the ring of invariants under the Borel subgroup consisting of the following 2 families: 
		\begin{enumerate}
			\item[(i)] $x_1^{q^m-1} x_2^{i (q-1)}$, $i \leq [m]_q$,
			\item[(ii)] $Q_{1,0}^{i_1} Q_{2,1}^{i_2}$, $i_1 < [m]_q$, $i_2 \leq [m-1]_q$. 
		\end{enumerate}

  It is a standard fact about the relative transfer $\tr \colon \Q_m(2)^B \to \Q_m(2)^G$ that it is onto since the index of $B$ in $G$ is relatively prime to $p$ (see \cite[Section 2.2]{Neusel-Smith-2002}), so a generating set for $\Q_m(2)^G$ can be obtained from the image of $\tr$. 
  
  In what follows, we adopt Goyal's notation $y_s$ for $\delta_2 (Q_{1,0}^s)$ for the sake of brevity. By direct inspection, we have the following identity in $\Q_m(2)$ for any $0 \leq s < [m]_q$:  
		\[
		x_1^{q^m-1} x_2^{s(q-1)} = y_0 x_1^{(s+1)(q-1)}  -  y_{s+1}.
		\]   
Since $y_s$ is already a $G_2$-invariant, this formula allows us to compute the transfer of the first family:   
		\[
		\tr (x_1^{q^m-1} x_2^{s(q-1)}) = y_0 \tr (x_1^{(s+1)(q-1)}) - y_{s+1}.
		\]
	Note that $\tr (x_1^{(s+1)(q-1)})$ is a Dickson polynomial divisible by $Q_{2,0}$ because $x_1^{(s+1)(q-1)}$ is a genuine $B$-invariant which vanishes when setting $x_1 = 0$. By the first two parts of Proposition \ref{prop: delta and Dickson}, the product of $y_0$ with $Q_{2,0}$ is trivial in $\Q_m(2)$. Thus 
    \[
    \tr (x_1^{q^m-1} x_2^{s(q-1)}) = -y_{s+1}.
    \]  
  For the second family, we argue similarly, noting that 
  \[
  \tr (Q_{1,0}^{i_1} Q_{2,1}^{i_2}) = \tr (Q_{1,0}^{i_1}) Q_{2,1}^{i_2} \in \D_2. 
  \]
The lemma follows. 
	\end{proof}

	The next step is to show that one only needs part of the Dickson algebra. The following observation about the Dickson monomials $Q_{2,1}^{\frac{q^{m-1}-q^i}{q-1}} Q_{2,0}^{\frac{q^i-1}{q-1}}$ at the "edge" of $\Delta^m_2$ is crucial.  
 \begin{proposition}\label{prop:y essential}
	For each $0 \leq i \leq  m-1$, we have the following decomposition in $S$:  
	\begin{equation*} 
		Q_{2,1}^{\frac{q^{m-1}-q^i}{q-1}} Q_{2,0}^{\frac{q^i-1}{q-1}} = \delta_{2} (Q_{1,0}^\frac{q^i-1}{q-1}) +   \text{essential monomials divisible by $Q_{2,0}^{\frac{q^{i+1}-1}{q-1}}$}. 
	\end{equation*}
\end{proposition}
\begin{proof}
	It is easy to see, by determinant manipulation, that $ \delta_{2} (Q_{1,0}^{\frac{q^i-1}{q-1}})  = y_{\frac{q^i-1}{q-1}}$ is a genuine Dickson polynomial. When $i=0$, the Dickson monomial $y_0$ must contain $Q_{2,1}^{\frac{q^{m-1}-1}{q-1}}$ as a nontrivial summand since both are reduced to $x_1^{q^m-q}$ when setting $x_2 = 0$. Furthermore, $y_0 - Q_{2,1}^{\frac{q^{m-1}-1}{q-1}}$ is divisible by $Q_{2,0}$. Note that in degree $q^m-q$,  $Q_{2,1}^{\frac{q^{m-1}-1}{q-1}}$ is the only non-essential Dickson monomial. We have proved that  
	\[
	y_0 = Q_{2,1}^{\frac{q^{m-1}-1}{q-1}} + \text{essential monomials divisible by $Q_{2,0}$}. 
	\]
For the general case, we make use of the identity $y_{m,\frac{q^i-1}{q-1}} =  Q_{2,0}^{\frac{q^i-1}{q-1}} y_{m-i,0}^{q^i}$ where $y_{m,s}$ stands for $\delta_{2;m}(Q_{1,0}^s)$. It follows that
	\[
	y_{\frac{q^i-1}{q-1}} =  Q_{2,1}^{\frac{q^{m-1}-q^i}{q-1}} Q_{2,0}^{\frac{q^i-1}{q-1}} + (fQ_{2,0})^{q^i}Q_{2,0}^{\frac{q^i-1}{q-1}} , 
	\]
	for some Dickson polynomial $f$. A Dickson monomial appearing in the second term of the sum above has the form $Q_{2,1}^{i_1} Q_{2,0}^{i_2}$ where $i_2 \geq \frac{q^{i+1}-1}{q-1}$. Comparing the degree, we have $i_1 \leq \frac{q^{m-1}-q^{i+1}}{q-1}$, hence $Q_{2,1}^{i_1} Q_{2,0}^{i_2}$ is essential. 
\end{proof}

\begin{remark}
It is easy to check that if $i > [m-1]_q$, then $Q_{2,0}^{i} = L_2^{i(q-1) \geq q^{m-1}} =0$ in $\Q_{m} (2)$. Also, by induction on $m$, we have another interesting formula: 
	\[
	Q_{2,1}^{\frac{2(q^{m-1}-1)}{q-1}} = y_{\frac{q^{m}-q}{q-1}} = x_1^{q^m-q} x_2^{q^m-q} \quad \text{in $\Q_{m}(2)$}. 
	\] 
\end{remark}

 We are now ready to finish the proof of Proposition \ref{prop:GL2}. 
 
 Consider a Dickson monomial $Q_{2,1}^{i_1} Q_{2,0}^{i_2}$ and suppose, by way of contradiction, that it is not in the span of $\mathcal{B}_{m}(2)$. Further, assume that it is the smallest monomial in grevlex order with this property. Since $Q_{2,0}^k = 0$ if $k > [m-1]_q$,  we must have $i_2 \leq \frac{q^{m-1}-1}{q-1} =[m-1]_q$. Let $i \leq m-1$ be the unique integer such that $[i]_q \leq i_2 < [i+1]_q$. Then $i_1 \geq \frac{q^{m-1}-q^i}{q-1}$. It follows from Proposition \ref{prop:y essential} that 
	\[
	Q_{2,1}^{i_1} Q_{2,0}^{i_2} - Q_{2,1}^{i_1 - \frac{q^{m-1}-q^i}{q-1}} Q_{2,0}^{i_2 - \frac{q^i-1}{q-1}} y_{\frac{q^i-1}{q-1}}
	\]
	is a sum of Dickson monomials in strictly smaller grevlex order than $(i_1, i_2)$, and at least one of them is again not in the span of $\mathcal{B}_{m}(2)$. We have a contradiction. 

\section{Proof of Proposition \ref{prop:G21}}\label{sec:G21} 
 It is straightforward to check that all elements in the six families listed in Proposition \ref{prop:G21} are $P(2,1)$-invariants. By results in Section \ref{sec: lower bound}, it suffices to show that $\mathcal{B}_m (2,1)$ is a generating set for $\Q_m(3)^{P(2,1)}$. 
 
 Suppose $f$ is a nonzero polynomial in $\Q_{m}(3)^{P(2,1)}$, and we write $f$ as a polynomial in $x_3$:
\[
f = x_3^{a(q-1)} f_{a}  + x_3^{(a-1)(q-1)} f_{a-1}   + \ldots,  
\]  
We will prove by downward induction on the highest $x_3$-degree $a(q-1)$.

Observe that the coefficients $f_i$ belong to the space of invariants $\Q_m(2)^G$ in the variables $x_1$ and $x_2$. If $f_a$ is a scalar multiple of $x_1^{q^m-1} x_2^{q^m-1}$, then 
since $x_3^{a(q-1)} x_1^{q^m-1} x_2^{q^m-1}$ is already $P(2,1)$-invariant (which is a member of the 5th and 6th family in $\mathcal{B}_m (2,1)$), we can subtract it from $f$ to obtain a new one with strictly smaller $x_3$-degree. 

Now let us assume that $f_a$, and hence all $f_i$, for degree reasons, are generated by $\delta_2(\Delta^m_1)$ and $\Delta^m_2$. In particular, its $x_1$-degree (and hence $x_2$-degree as well, by symmetry) is strictly less than $q^m-1$. We claim that $a$ must be a multiple of $q$. Indeed, since $f$ is a $B$-invariant, from the list of generators given in  \ref{subsec: Borel}, we see that $f$ must be a linear combination of polynomials of the form $x_1^{j_1 (q-1)} \delta_2 (Q_{2,1}^{j_2})$ or $x_1^{j_1 (q-1)} Q_{2,1}^{j_2} Q_{3,2}^{j_3}$. In both cases, the highest $x_3$-degrees are divisible by $q$. 

If $f_a$ contains $\delta_2 (Q_{1,0}^{i_1})$ for some $i_1 <[m]_q$, then we can subtract from $f$ an appropriate polynomial from the third or fourth family in the list \ref{prop:G21}. So we can assume that $f_a \in \Delta^m_2$. In this case, it follows that $x_3^{a(q-1} f_a$ is 
a nontrivial summand of a $B$-invariant polynomial of the form $x_1^{j_1 (q-1)} Q_{2,1}^{j_2} Q_{3,2}^{j_3}$. In particular, $a$ must be divisible by $q^2$. 

Finally, when $a(q-1) = iq^2(q-1)$ for some $i< [m-2]_q$ or $a(q-1) = q^m-q^2$, by subtracting from $f$ an appropriate linear combination of the first two families of $\mathcal{B}_{m}(2,1)$, we will obtain a $P(2,1)$-invariant polynomial of strictly smaller $x_3$-degree. The result follows by induction on the highest $x_3$-degree. \qed 
 
\begin{remark}
We can also prove that $\mathcal{B}_m (2,1)$ is linearly independent as follows. Consider these elements of $\mathcal{B}_m (2,1)$ as polynomials in $x_3$. The table below records the coefficient of the highest $x_3$-degree. 
				\begin{center} 
					\begin{tabular}{|l|c|c|}
						\hline
						$\mathcal{B}_{m}(2,1)$	& Highest $x_3$-degree & Coefficient  \\
						\hline
						(1) $Q_{2,1}^{i_1} Q_{2,0}^{i_2}Q_{3,2}^i $	& $i(q^3-q^2) < q^m-q^2$ & $Q_{2,1}^{i_1} Q_{2,0}^{i_2}$,  $(i_1, i_2) \in \Delta^{m}_2$   \\
						\hline
						(2) $Q_{2,1}^{i_1} Q_{2,0}^{i_2}\delta_3 (1) $	& $ q^m-q^2$ & $Q_{2,1}^{i_1} Q_{2,0}^{i_2}$,  $(i_1, i_2) \in \Delta^m_2$   \\
						\hline
						(3)	$\delta_2 (Q_{1,0}^{i_1} Q_{2,1}^{i})$ & $i (q^2-q) < q^m-q$ &  $\delta_2(Q_{1,0}^{i_1})$, $i_1<[m]_q$\\
						\hline
						(4)	$\delta_2 (Q_{1,0}^{i_1} \delta_2(1))$ & $q^m-q$ & $\delta_2(Q_{1,0}^{i_1})$, $i_1<[m]_q$ \\
						\hline
						(5)	$\delta_1(\delta_1(Q_{1,0}^i))$ & $i(q-1)<q^m-1$ & $x_1^{q^m-1} x_2^{q^m-1}$ \\
						\hline
						(6)	$\delta_1(\delta_1(\delta_1(1)))$ & $q^m-1$ & $x_1^{q^m-1} x_2^{q^m-1}$  \\
						\hline
					\end{tabular}
				\end{center} 
From the table, it is easy to see that the polynomials in $\mathcal{B}_m (2,1)$ are linearly independent. 		 
\end{remark}			
			\section{Proof of Proposition \ref{prop:G12}}\label{sec:G12}
It suffices to show that $\mathcal{B}_m(1,2)$ is a generating set, and the proof is by downward induction on the lowest $x_1$-degree. Write 
\[
f = x_1^{i(q-1)} g_i (x_2, x_3) + x_1^{(i+1)(q-1)} g_{i+1} (x_2, x_3) + \ldots 
\]
Note that all coefficient polynomials $g_j$ are $G_2$-invariant in $x_2$ and $x_3$. If $i(q-1) = q^m-1$, then $x_1^{q^m-1} g(x_2, x_3)$ is a linear combination of polynomials from the last three families of $\mathcal{B}_m (1,2)$. We can then subtract from $f$ this polynomial and can now assume that $i(q-1)  \leq q^m-q$.  

Since $f$ is a $B$-invariant, it is a linear combination of polynomials of the form $x_1^{j_1(q-1)} \delta_2 (Q_{2,1}^{j_2})$ or $Q_{1,0}^{j_1} Q_{2,1}^{j_2} Q_{3,2}^{j_3}$. Their lowest $x_1$-degree summands are $x_1^{j_1 (q-1)} x_2^{q^m-q} x_3^{j_2q(q-1)}$ and $x_1^{j_1 (q-1)} x_2^{j_2 q(q-1)} Q_{2,1}^{qj_3} (x_2, x_3)$ respectively. In both cases, the exponents of $x_2$ and $x_3$ in these summands are divisible by $q(q-1)$. So the same is true for $g_i (x_2, x_3)$. Since $g_i \in \Q_{m}(x_2, x_3)^G$, it must be a $q$-th power of a polynomial in $\Q_{m-1}(2)$ in the variables $x_2, x_3$. 

 On the other hand, the lowest $x_1$-degree summand of the first three families of $\mathcal{B}_m(1,2)$ are 
 \begin{enumerate}
     \item[(1)] $x_1^{j_1(q-1)} \big(Q_{2,1}^{i_1} (x_2,x_3) Q_{2,0}^{i_2}(x_2,x_3)\big)^q$, $(i_1, i_2) \in \Delta^{m-1}_2$;
     \item[(2)] $x_1^{j_1 (q-1)} \big(\delta_{2;m-1} (Q_{1,0}^{j_2})\big)^q$, $j_2 < [m-1]_q$;
     \item[(3)] $x_1^{j_1 (q-1)} \big(x_2^{q^{m-1}-1} x_3^{q^{m-1}-1}\big)^q$. 
 \end{enumerate}
Hence $x_1^{i(q-1)} g_i (x_2, x_3)$ is the lowest $x_1$-degree part of a polynomial in the span of $\mathcal{B}_m (1,2)$. By subtracting from $f$ an appropriate linear combination of polynomials in $\mathcal{B}_m (1,2)$, we obtain another $P(1,2)$-invariant polynomial whose lowest $x_1$-degree is strictly greater than that of $f$. By induction, the required result follows. \qed

		\begin{remark}\label{rem:G12}
We can also prove directly that $\mathcal{B}_m(1,2)$ is linearly independent. In fact, the generators can be distinguished by looking at the lowest $x_1$-degrees and the corresponding coefficients as given in the table below. 
				\begin{center} 
					\begin{tabular}{|l|c|c|}
						\hline
						$\mathcal{B}_m(1,2)$	& Lowest $x_1$-degree & Coefficient  \\
						\hline
						(1) $Q_{1,0}^{j_1}  Q_{3,2}^{i_1} Q_{3,1}^{i_2}$ & $j_1(q-1)< q^m-1$  & $Q_{2,1}^{qi_1}(x_2, x_3) Q_{2,0}^{qi_2} (x_2, x_3)$, $(i_1, i_2) \in \Delta^{m-1}_2$ \\
						\hline 
						(2)	$Q_{1,0}^{j_1} \delta_3 (Q_{2,1}^{j_2})$ & $j_1 (q-1) < q^m-1$  & $y_{m-1,j_2}^q (x_2, x_3)$, $j_2<[m-1]_q$ \\
						\hline 
						(3)	$Q_{1,0}^{j_1} \delta_2(\delta_2(1))$& $j_1(q-1) < q^m-1$ & $x_2^{q^m-q} x_3^{q^m-q}$\\
						\hline 
						(4) $\delta_1(Q_{2,1}^{i_1}Q_{2,0}^{i_2})$	& $q^m-1$ & $Q_{2,1}^{i_1}(x_2, x_3) Q_{2,0}^{i_2} (x_2, x_3)$, $(i_1, i_2) \in \Delta^m_2$ \\
						\hline 
						(5) $\delta_1(\delta_2(Q_{1,0}^{j_1}))$	& $q^m-1$ &  $y_{j_1} (x_2, x_3), j_1<[m]_q$ \\
						\hline
						(6)	$\delta_1(\delta_1(\delta_1(1)))$ & $q^m-1$ & $x_2^{q^m-1} x_3^{q^m-1}$\\
						\hline
					\end{tabular}
				\end{center} 
			 	
   \end{remark} 	
			\section{Proof of Proposition \ref{prop:GL3}}\label{sec:GL3}
We need to prove that $\mathcal B_m (3) $ is a generating set for $\Q_m (3)^G $. Our strategy is the same as in the rank $2$ case, but the details are much more involved. The proof is divided into three steps. Using the transfer, we first construct a slightly bigger generating set $\mathcal{B}'$ of the required form without restricting the Dickson polynomials:
\[
\mathcal{B}' = \delta_1^3 (\Delta^m_0) \coprod  \delta_2^2 (\Delta^m_1) \coprod  \delta_3 (\D_2) \coprod  \D_3.  
\]
Once this is done, we then show that $\delta_3 (\D_2)$ is actually contained in the span of 
\[
\delta_1^3 (\Delta^m_0) \coprod  \delta_2^2 (\Delta^m_1) \coprod  \delta_3 (\Delta^m_2),
\]
and as a result, $\D_2$ can be replaced by the much smaller subset $\Delta^m_2$. Finally, we prove that when restricted to $\Q_m(3)$, the rank $3$ Dickson algebra belongs to the subspace spanned by $\mathcal{B}_m(3)$. This is proved by first showing that all the edge monomials of $\Delta^m_3$ are in this span. Thus the subspace spanned by $\mathcal{B}_m(3)$ is a $\D_3$-submodule of $\Q_m(3)^G$. Since it contains $\Delta^m_3$, it will contain the entirety of $\D_3$.

 \begin{lemma}\label{lem: bigger gen set} The set $\mathcal{B}'$ consisting of 4 families below is a generating set for $\Q_m(3)^G$.
\begin{enumerate}
					\item $x_1^{q^m-1} x_2^{q^m-1} x_3^{q^m-1}$. 
					\item $\delta_2^2 (\Delta^m_1) = \{ a_{m,3,s}$, $0 \leq s < [m]_q\}$. 
					\item $\delta_3 (\D_2)$.  
					\item $\D_3$. 
				\end{enumerate}
 \end{lemma}
\begin{proof}
A generating set for $\Q_m(3)^G$ can be obtained by taking the image of the relative transfer $$\tr_{P(2,1)}^G \colon  \Q_m^{P(1,2)} \longrightarrow \Q_m^G$$ from the parabolic subgroup $P(1,2)$ to $G$, applied to the basis of $\Q_m^{P(1,2)}$ given in Proposition \ref{prop:G12}. 

First, observe that $Q_{1,0}^{j_1} Q_{3,2}^{i_1} Q_{3,1}^{i_2}$ is a genuine $P(1,2)$-invariant, and hence its transfer is a genuine $G$-invariant, i.e., a Dickson polynomial. 

Next, consider the second family. Since $\delta_3 (Q_{2,1}^{j_2})$ is already a $G$-invariant in $\Q_m(3)$, we have 
\[
\tr_{P(2,1)}^G (Q_{1,0}^{j_1} \delta_3 (Q_{2,1}^{j_2})) = (\tr_{P(2,1)}^G Q_{1,0}^{j_1}) \delta_3 (Q_{2,1}^{j_2}). 
\]
On the other hand, the transfer $\tr_{P(2,1)}^G (Q_{1,0}^{j_1})$ is a genuine Dickson polynomial divisible by $Q_{3,0}$ unless $j_1 = 0$. Therefore, if $j_1$ is positive, the whole expression vanishes in $\Q_m(3)$:  
 \[
 \tr_{P(2,1)}^G (Q_{1,0}^{j_1} \delta_3 (Q_{2,1}^{j_2})) = 0 \quad \text{in $\Q_m$}. 
 \]
A similar argument applies to the third family. Since $\delta_2^2 (1)$ is already a $G$-invariant and is annihilated by $Q_{3,0}$, its transfer also vanishes. Hence from the second and third families in $\mathcal{B}_m (1,2)$, only elements  $\delta_3 (Q_{2,1}^{j_2})$ (for $j_2 <[m-1]_q$) and $\delta_2^2 (1)$ contribute nontrivially.  

Now we turn to the fourth family. Recall that 
\[
\delta_1 (Q_{2,1}^{i_1} Q_{2,0}^{i_2})= x_1^{q^m-1} Q_{2,1}^{i_1}(x_2, x_3) Q_{2,0}^{i_2} (x_2, x_3). 
\]
We claim that 
\[
\tr_{P(1,2)}^G  (\delta_1 (Q_{2,1}^{i_1} Q_{2,0}^{i_2})) = \delta_3 (Q_{2,1}^{i_1} Q_{2,0}^{i_2+1}). 
\]
Indeed, by Proposition~\ref{prop:delta3}, the right-hand side is $G$-invariant. Moreover, we use the identity: 
\[
V_3 (x_2,x_1,x_3) = V_2 (V_{2}(x_2,x_1), V_2 (x_2,x_3)) = V_2(x_2,x_3)^q - V_2^{q-1}(x_2,x_1) V_2(x_2,x_3), 
\]
from which it follows that 
\[
L_2^q (x_2,x_3) = x_2^{q} V_2(x_2,x_3)^q = x_2^q V_3 (x_2,x_1,x_3) +  Q_{2,0}(x_1,x_2) L_2(x_2,x_3). 
\]
Applying the Laplace expansion along the last row of the numerator, we obtain: 
\[
\delta_3 (Q_{2,1}^{i_1} Q_{2,0}^{i_2+1}) = Q_{2,0} \delta_3 (Q_{2,1}^{i_1} Q_{2,0}^{i_2}) - \delta_2 (Q_{1,0} Q_{2,1}^{i_1} Q_{2,0}^{i_2}). 
\]
Now consider $x_1^{q^m-1} Q_{2,1}^{i_1}(x_2, x_3) Q_{2,0}^{i_2} (x_2, x_3)$ as a $B_3$-invariant and take the transfer from $B_3$ to $P(2,1)$. From Section~\ref{sec:G2}, we have
						\[
						\tr_{B_2}^{G_2} (x_1^{q^m-1} x_2^{a(q-1)}) = - y_{a+1} (x_1, x_2) = - \delta_2 (Q_{1,0}^{a+1}).
						\]
 Write 
\[
Q_{2,1}^{i_1} (x_2, x_3) Q_{2,0}^{i_2} (x_2, x_3) = \sum_{(a,b)} \lambda (a,b) x_2^{a(q-1)} x_3^{b(q-1)}, \quad \lambda(a,b) \in \mathbb{F}_q,   
\]
then 
\begin{align*}
\tr_{B_3}^{P(2,1)} \big(x_1^{q^m-1} Q_{2,1}^{i_1} (x_2, x_3) Q_{2,0}^{i_2} (x_2, x_3) \big) =&    - \sum_{(a,b)} \lambda(a,b) \delta_2 (Q_{1,0}^{a+1}) x_3^{b(q-1)}\\
=&  -   \delta_2 \big(\sum_{(a,b)} \lambda(a,b) x_1^{(a+1)(q-1)} x_2^{b(q-1)}\big )\\
=& - \delta_2 (Q_{1,0} Q_{2,1}^{i_1} Q_{2,0}^{i_2})\\
=& \delta_3 (Q_{2,1}^{i_1} Q_{2,0}^{i_2+1})-Q_{2,0} \delta_3 (Q_{2,1}^{i_1} Q_{2,0}^{i_2}). 
\end{align*}
Finally, we observe that 
\[
\tr_{P(2,1)}^G \big(Q_{2,0} \delta_3 (Q_{2,1}^{i_1} Q_{2,0}^{i_2})\big) = \tr_{P(2,1)}^G (Q_{2,0}) \delta_3 (Q_{2,1}^{i_1} Q_{2,0}^{i_2}) = 0, 
\]
which proves the claim for the fourth family. 

Now we consider the fifth family (note that the 6th family already consists of $G$-invariants). A similar argument shows that if $j < [m]_q$, then 
				\[
				\tr_{P(1,2)}^G  (x_1^{q^m-1} \delta_2(Q_{1,0}^{j})(x_2,x_3)) = -\delta_2^2 (Q_{1,0}^{j+1}).   
				\]
    Indeed, recall 
    \[
						\tr_{B_3}^{P(2,1)} (x_1^{q^m-1} y_{j} (x_2,x_3)) =    - \sum_{a=j}^{(q^m-q)/(q-1)} y_{a+1} (x_1, x_2) x_3^{q^m-q + j(a-1) - a(q-1)}. 
						\]
   This sum is exactly $- a_{m,3,j+1}=- \delta_2^2(Q_{1,0}^{j+1})$, which is a $G_3$-invariant. The result follows by the transitivity of the transfer. The proof is complete. 		
\end{proof}

The next step is to show that in the third family, one does not need to use all $\D_2$. The following technical lemma will be used repeatedly. 
\begin{lemma}\label{lem: delta2 f delta2}
    Let $s, t, i$ be non-negative integers. Then the following equality holds in $\Q_m(3)$:
\[
\delta_{2} \big(V_1^{s(q-1)} V_2^{t(q-1)} \delta_{2} (Q_{1,0}^i)\big) = 
\begin{cases}
0 & \text{if $t \geq 1$ and $s > 1$};\\
x_1^{q^m-1} x_2^{q^m-1} x_3^{(qt+i-1)(q-1)} & \text{if $t \geq 1$ and $s =1$};\\
(t+1) x_1^{q^m-1} x_2^{q^m-1} x_3^{(qt+i-2)(q-1)} & \text{if $t \geq 1$ and $s =0$};\\
\delta_{2}^2 (Q_{1,0}^{i+s}) & \text{if $t=0, s \leq 1$};\\
\delta_{2}^2 (Q_{1,0}^{i+s}) -x_1^{q^m-1} x_2^{q^m-1} x_3^{(s+i-2)(q-1)} & \text{if $t=0, s \geq 2$}.
\end{cases}
\]
\end{lemma}
\begin{proof}
Recall that $\delta_{2} (f)$ is a polynomial if $f$ is a polynomial (possibly with more than one variable), which is $G_1$-invariant. Hence, the expression under consideration is a polynomial.  
If $t \geq 1$, unravelling the definition, we see that $\delta_{2} (V_1^{s(q-1)} V_2^{t(q-1)} \delta_{2} (Q_{1,0}^i))$ equals 
\begin{multline*}
\frac{1}{L_2 (x_1,x_2)} \bigg[x_2^{q^m} x_1^{s(q-1)} V_2^{t(q-1)-1} (x_1,x_3) \begin{vmatrix}
    x_1 & x_3\\
    x_1^{q^m} x_3^{i(q-1)} & x_3^{q^m} x_1^{i(q-1)} 
\end{vmatrix} -\\ x_1^{q^m} x_2^{s(q-1)} V_2^{t(q-1)-1} (x_2,x_3) \begin{vmatrix}
    x_2 & x_3\\
    x_2^{q^m} x_3^{i(q-1)} & x_3^{q^m} x_2^{i(q-1)} 
\end{vmatrix}\bigg]\\
= x_3^{q^m} \; \frac{x_2^{q^m} x_1^{(s+i)(q-1)+1} V_2^{t(q-1)-1} (x_1,x_3)-x_1^{q^m} x_2^{(s+i)(q-1)+1} V_2^{t(q-1)-1} (x_2,x_3)}{L_2 (x_1,x_2)} -\\  x_1^{q^m} x_2^{q^m} x_3^{i(q-1)+1} \; \frac{x_1^{s(q-1)} V_2^{t(q-1)-1}(x_1, x_3) - x_2^{s(q-1)} V_2^{t(q-1)-1}(x_2, x_3)}{L_2(x_1,x_2)}.
\end{multline*}
The first summand is trivial in $\Q_m(3)$ because the quotient involved is a polynomial. Indeed, in this quotient, the numerator is divisible by $x_1$ and vanishes whenever there is a linear relation $x_2 = a_1 x_1$. 

By a similar argument, the second summand is also trivial if $s >1$, and equals $x_1^{q^m-1} x_2^{q^m-1} x_3^{(qt+i-1)(q-1)}$ if $s = 1$. If $s=0$, it equals
\[
(t(q-1)-1) x_1^{q^m-1} x_2^{q^m-1} x_3^{(qt+i-2)(q-1)}. 
\]

Now suppose $t=0$, we get 
 \begin{multline*}
 \frac{1}{L_2 (x_1, x_2)} \big[x_1^{q^m}  x_2^{s(q-1)}  \frac{\begin{vmatrix}
		x_2 & x_3\\ x_2^{q^m} x_3^{i(q-1)} & x_3^{q^m} x_2^{i(q-1)}  
\end{vmatrix}}{V_2 (x_2,x_3)} -x_2^{q^m}  x_1^{s(q-1)}   \frac{\begin{vmatrix}
		x_1  & x_3 \\ x_1^{q^m}x_3^{i(q-1)} & x_3^{q^m} x_1^{i(q-1)}
\end{vmatrix}}{V_2 (x_1,x_3)}\big]\\
= \frac{1}{L_2 (x_1, x_2)} \big[x_1^{q^m} x_2  \frac{\begin{vmatrix}
		x_2  & x_3 \\ x_2^{q^m} x_3^{(s+i)(q-1)}  & x_3^{q^m}  x_2^{(s+i)(q-1)} 
\end{vmatrix}}{L_2 (x_2,x_3)} -x_2^{q^m}  x_1   \frac{\begin{vmatrix}
		x_1 & x_3\\ x_1^{q^m} x_3^{(s+i)(q-1)}  & x_3^{q^m}  x_1^{(s+i)(q-1)} 
\end{vmatrix}}{L_2 (x_1,x_3)}\big] +\\
\frac{x_1^{q^m}x_2^{q^m}x_3^{i(q-1)}}{L_2 (x_1,x_2)} \Bigg[\frac{x_2^{s(q-1)}-x_3^{s(q-1)}}{x_2^{q-1}-x_3^{q-1}} -  \frac{x_1^{s(q-1)}-x_3^{s(q-1)}}{x_1^{q-1}-x_3^{q-1}}  \Bigg].
\end{multline*} 
The first summand is $\delta_2^2 (Q_{1,0}^{s+i})$. The second summand is trivial when $s\leq 1$ and equals $-x_1^{q^m-1} x_2^{q^m-1} x_3^{(s-2+i)(q-1)}$ if $s \geq 2$. 
\end{proof}
\begin{corollary}\label{cor 73}
If $g \in \D_2$ and $i \geq 0$, then $\delta_3 (g \delta_2 (Q_{1,0}^i))$ is a linear combination of polynomials from the following families: 
\begin{itemize}
    \item $\delta_2^{2} (Q_{1,0}^s)$, $0 \leq s <[m]_q$. 
    \item $x_1^{q^m-1} x_2^{q^m-1} x_3^{s(q-1)}$, $0 \leq s \leq [m]_q$. 
\end{itemize}
\end{corollary}
 \begin{proof}
   We first observe that by direct inspection, $\delta_{3} (g \delta_{2} (Q_{1,0}^i))$ can be written as a quotient of determinants: 
		\[
	\delta_{3} (g \delta_{2} (Q_{1,0}^i)) = 	\frac{ 
		\begin{vmatrix}
			x_1^{i(q-1)+1} & x_2^{i(q-1)+1} & x_3^{i(q-1)+1}\\
			x_1^{q^m} & x_2^{q^m} & x_3^{q^m}\\
			x_1^{q^m} g(x_2, x_3) &  x_2^{q^m} g(x_1, x_3) & x_3^{q^m} g(x_1, x_2) 
		\end{vmatrix} 
	}{L_3 (x_1, x_2, x_3)}.
		\]			
Since
 \begin{align*}
 L_3 (x_1, x_2, x_3) =& V_3 (x_1, x_2, x_3) L_2 (x_1, x_2)\\
 =& [V_2 (x_2, x_3)^{q-1}-V_2 (x_2, x_1)^{q-1}] V_2 (x_2,x_3) L_2 (x_1, x_2)\\
 =& [V_2 (x_1, x_3)^{q-1}-V_2 (x_1, x_2)^{q-1}] V_2 (x_1,x_3) L_2 (x_1, x_2),
 \end{align*}
 the quotient above can be rewritten in the form: 
\begin{align*} 
 \varphi (h) :=& \frac{1}{L_2 (x_1, x_2)} \bigg[x_1^{q^m} x_2 h(x_1,x_2,x_3) \frac{\begin{vmatrix}
 		x_2^{i(q-1)+1} & x_3^{i(q-1)+1}\\ x_2^{q^m} & x_3^{q^m} 
 \end{vmatrix}}{L_2 (x_2,x_3)} - x_2^{q^m} x_1 h(x_2, x_1, x_3) \frac{\begin{vmatrix}
 x_1^{i(q-1)+1} & x_3^{i(q-1)+1}\\ x_1^{q^m} & x_3^{q^m} 
\end{vmatrix}}{L_2 (x_1,x_3)}\bigg]\\
= & \frac{1}{L_2 (x_1, x_2)} (x_1^{q^m} x_2 h(x_1,x_2,x_3) \delta_2 (Q_{1,0}^i)(x_2,x_3) - x_2^{q^m} x_1 h(x_2, x_1, x_3) \delta_2 (Q_{1,0}^i)(x_1,x_3)),
\end{align*}
where
\[
h(x_1, x_2, x_3) = \frac{g(x_2, x_3)-g(x_2,x_1)}{V_2 (x_2, x_3)^{q-1} -V_2 (x_2, x_1)^{q-1}}. 
\]
Since $g \in \D_2$, it can be written in terms of upper triangular invariants 
$$g = \sum_{s,t} V_1^{s(q-1)} V_2^{t(q-1)},$$ where the sum is taken over certain subset $A \subset \mathbb{N}^2$. Then 
\[
h(x_1, x_2, x_3) = \sum_{(s,t) \in A, t \geq 1} x_2^{s(q-1)} \bigg(\sum_{j=0}^{t-1} V_2(x_2,x_3)^{(t-1-j)(q-1)} V_2(x_2,x_1)^{j(q-1)}\bigg).  
\]
Observe that if $j> 0$, then both $x_1^{q^m}x_2V_2(x_2,x_1)^{j(q-1)}/L_2$ and $x_2^{q^m}x_1V_2(x_1,x_2)^{j(q-1)}/L_2$ are zero in $\Q_m$, so the terms in $h$ corresponding to those $j$ do not contribute to $\varphi(h)$.
It follows that $\varphi(h)$ is reduced to the sum 
\[
\sum_{(s,t) \in A, t \geq 1} \varphi (x_2^{s(q-1)} V_2 (x_2, x_3)^{(t-1)(q-1)}). 
\]
On the other hand, we have an equality which is essentially by definition: 
\[
\varphi (x_2^{s(q-1)} V_2 (x_2, x_3)^{(t-1)(q-1)}) = \delta_{2} (V_1^{s(q-1)} V_2^{(t-1)(q-1)} \delta_{2} (Q_{1,0}^{i})).
\]
We can now use the computation in Lemma \ref{lem: delta2 f delta2} of the righ-hand side to finish the proof.  
 \end{proof}

We are now ready to replace $\D_2$ by $\Delta^m_2$ in the generating set $\mathcal B'$ of $\Q_m(3)^G$. 
\begin{lemma}\label{lem:Dickson submodule delta3}
The subspace of $\Q_m(3)^G$ spanned by the following three families of elements 
\begin{enumerate}
    \item $\delta_1^3 (\Delta^m_0)$, 
    \item $\delta_2^2 (\Delta^m_1)$,
    \item $\delta_3 (\Delta^m_2)$, 			  
\end{enumerate}
contains $\delta_3 (\D_2)$.  
\end{lemma}
\begin{proof}
Suppose $f$ is a Dickson monomial in $\D_2$ but is not essential. Then it can be written in the form 
\begin{equation*}
Q_{2,1}^{\frac{q^{m-1}-q^i}{q-1}} Q_{2,0}^{\frac{q^i-1}{q-1}} Q_{2,1}^a Q_{2,0}^b,
\end{equation*}
for some $0 \leq i \leq m-1$, $0 \leq a,b$.  
 Recall from \ref{prop:y essential} that 
\[
Q_{2,1}^{\frac{q^{m-1}-q^i}{q-1}} Q_{2,0}^{\frac{q^i-1}{q-1}} = \delta_{2} (Q_{1,0}^{\frac{q^i-1}{q-1}}) + \text{essential monomials divisible by $Q_{2,0}^{\frac{q^{i+1}-1}{q-1}}$}. 
\]
Repeating this process if necessary, we can write $f$ as a sum of essential monomials and Dickson polynomials of the form
\[
\delta_{2} (Q_{1,0}^{\frac{q^j-1}{q-1}}) Q_{2,1}^{a_j} Q_{2,0}^{b_j}, \quad j \geq i. 
\]
Thus it suffices to consider $\delta_3(f)$ for $f = g \delta_{2} (Q_{1,0}^{\frac{q^i-1}{q-1}})$, where $g \in \D_2$ and $0 \leq i \leq m-1$.

We know from Corollary \ref{cor 73} that $\delta_3 (g \delta_2 (Q_{1,0}^{\frac{q^i-1}{q-1}}))$ is 
a linear combination of polynomials in the following families:  \begin{itemize}
    \item $\delta_2^{2} (\Delta^m_1)$,
    \item $x_1^{q^m-1} x_2^{q^m-1} x_3^{s(q-1)}$, $0 \leq s \leq [m]_q$. 
\end{itemize} 
 Since $\delta_2 (Q_{1,0}^{\frac{q^i-1}{q-1}})$ is a genuine Dickson polynomial, $\delta_{3} (g \delta_{2} (Q_{1,0}^{\frac{q^i-1}{q-1}}))$ is a $G_3$-invariant in $\Q_m(3)$. We conclude that $\delta_3 (g \delta_2 (Q_{1,0}^{\frac{q^i-1}{q-1}}))$ must belong to the space
\[
 \Span \{ \delta_{2}^2 (\Delta^m_1), \delta_{1}^3 (1)\},
\]
because any monomial of the form $x_1^{q^m-1} x_2^{q^m-1} x_3^{i(q-1)}$, $0 \leq i< [m]_q$, is not $G_3$-invariant in $\Q_m(3)$.  
 \end{proof}

 The third, and final step is to replace the full Dickson algebra $\D_3$ by $\Delta^m_3$. We first need some preliminary calculations. 
 In the following, $\delta_{3}$ will again stand for $\delta_{3;m}$.
 
					 				
\begin{lemma}\label{lem: delta3 formula} The following equalities hold in $S$:   
						\begin{enumerate}
				\item 	$\delta_{3} (Q_{2,0}^{\frac{q^{i}-1}{q-1}}) = \big(\delta_{3;m-i}(1)\big)^{q^i} Q_{3,0}^{\frac{q^{i}-1}{q-1}}$ for all $i\ge 0$. In particular, 
    $$\delta_{3} (Q_{2,0}^{\frac{q^{m-2}-1}{q-1}}) =Q_{3,0}^{\frac{q^{m-2}-1}{q-1}}\text{\quad and \quad }\delta_{3}(Q_{2,0}^{\frac{q^{m-3}-1}{q-1}}) = Q_{3,2}^{q^{m-3}} Q_{3,0}^{\frac{q^{m-3}-1}{q-1}}.$$
    \item  $Q_{3,0}^k = 0$ in $\Q_m$ for all $k >\frac{q^{m-2}-1}{q-1}$. 
    \item 	$\delta_{3}(Q_{2,1}^{q^{m-3}} Q_{2,0}^{\frac{q^{m-3}-1}{q-1}}) = Q_{3,1}^{q^{m-3}}  Q_{3,0}^{\frac{q^{m-3}-1}{q-1}}$.
				\item  In general, for  $\lambda_2 \geq \lambda_3 \geq 0$, then   
				\[
				\delta_{3} \big(Q_{2,1}^{\frac{q^{\lambda_2}-q^{\lambda_3}}{q-1}} Q_{2,0}^{\frac{q^{\lambda_3}-1}{q-1}} \big) = \bigg(\delta_{3;m-\lambda_3} \big( Q_{2,1}^{\frac{q^{\lambda_2-\lambda_3}-1}{q-1} } \big)\bigg)^{q^{\lambda_3}}  Q_{3,0}^{\frac{q^{\lambda_3}-1}{q-1}}. 
				\] 
						\end{enumerate}
					\end{lemma}
     \begin{proof}
    The proofs of these statements are quite similar and essentially make use of the Laplace expansion. We will prove the last statement only. Apply the Laplace expansion along the last row of the determinant in the numerator of $\delta_{3} \big(Q_{2,1}^{\frac{q^{\lambda_2}-q^{\lambda_3}}{q-1}} Q_{2,0}^{\frac{q^{\lambda_3}-1}{q-1}} \big)$, we obtain 
\[
x_1^{q^m} Q_{2,1}^{\frac{q^{\lambda_2}-q^{\lambda_3}}{q-1}} (x_2,x_3) L_2^{q^{\lambda_3}} (x_2,x_3)  - x_2^{q^m} Q_{2,1}^{\frac{q^{\lambda_2}-q^{\lambda_3}}{q-1}} (x_1, x_3) L_2^{q^{\lambda_3}} (x_1,x_3)  + x_3^{q^m} Q_{2,1}^{\frac{q^{\lambda_2}-q^{\lambda_3}}{q-1}} (x_1,x_2) L_2^{q^{\lambda_3}} (x_1,x_3). 
\]
After writing this expression as a $q^{\lambda_3}$ power, the result follows easily.  
     \end{proof}
Our second technical result is similar to Proposition \ref{prop:y essential} in the case of 3 variables. 

\begin{proposition}\label{prop:key reduction 3} 
Let $m \geq 2$ be an integer. 
 \begin{enumerate}
     \item[(i)] The following equality holds in $S$:
 \begin{equation*}
 \frac{[0,1,m-1]}{[0,1,2]} \; \delta_{3;m} (Q_{2,1}) - \frac{[0,2,m-1]}{[0,1,2]} \; \delta_{3;m} (1) = Q_{3,1}^{\frac{q^{m-2}-1}{q-1}} + \text{other terms},
 \end{equation*} 
 where the other terms are in the ideal $(Q_{3,0})$ of $\D_3$ generated by $Q_{3,0}$. 
 \item[(ii)] For $0 \leq \ell \leq m-2$, the following equality holds in $S$:
 \[
 \frac{[0,1,\ell + 1]}{[0,1,2]} \; \delta_{3;m} (Q_{2,1}) - \frac{[0,2,\ell +1]}{[0,1,2]} \; \delta_{3;m} (1) = Q_{3,2}^{\frac{q^{m-2}-q^{\ell}}{q-1}} Q_{3,1}^{\frac{q^{\ell}-1}{q-1}} + \text{other terms},
 \]  
 where the other terms are in the ideal $(Q_{3,1}^{\frac{q^{\ell+1}-1}{q-1}},Q_{3,0})$ of $\D_3$. 
 \end{enumerate}
\end{proposition}  
 	\begin{proof}
 The key observation is that if two Dickson polynomials coincide after setting one of the variables, say $x_3$, to zero, then their difference must be a multiple of $Q_{3,0}$.  

Note that in Dickson's notation, we can write 
\[
\delta_{3;m} (Q_{2,1}) = \frac{[0,2,m]}{[0,1,2]}, \quad \delta_{3;m}(1) = \frac{[0,1,m]}{[0,1,2]}, 
\]
so that all the terms in the above equations are genuine Dickson polynomials. 

\noindent\textbf{Proof of (i):} 
When $x_3 =0$, the equation under consideration becomes 
\[
\frac{[1,m-1]}{[1,2]} \cdot \frac{[2,m]}{[1,2]} - \frac{[2,m-1]}{[1,2]} \cdot \frac{[1,m]}{[1,2]} = Q_{2,0}^{\frac{q^{m-1}-q}{q-1}}. 
\]
It is now a simple computation that the two sides are the same. 




\noindent \textbf{Proof of (ii):} Again, we consider the reduction of the left-hand side when setting $x_3=0$. The result, after simplification, is 
\[
\delta_{2;m-\ell-1}(1)^{q^{\ell+1}} Q_{2,0}^{\frac{q^{\ell+1}-q}{q-1}}. 
\]
Now setting $x_2=0$, $\delta_{2;m-\ell-1}(1)$ is reduced to $Q_{1,0}^{\frac{q^{m-\ell-1}-q}{q-1}}$, which implies that
\[
\delta_{2;m-\ell-1}(1) = Q_{2,1}^{\frac{q^{m-\ell-2}-1}{q-1}} + \text{other terms divisible by $Q_{2,0}$}.  
\]
It follows that 
\[
\delta_{2;m-\ell-1}(1)^{q^{\ell+1}} Q_{2,0}^{\frac{q^{\ell+1}-q}{q-1}} = Q_{2,1}^{\frac{q^{m-1}-q^{\ell+1}}{q-1}} Q_{2,0}^{\frac{q^{\ell+1}-q}{q-1}} + \text{other terms divisible by $Q_{2,0}^{\frac{q^{\ell+2}-q}{q-1}}$.}
\]
Thus, we can conclude that the left-hand side equals 
\[
Q_{3,2}^{\frac{q^{m-2}-q^{\ell}}{q-1}} Q_{3,1}^{\frac{q^{\ell}-1}{q-1}} + \text{other terms},
 \]  
 where the other terms consist of Dickson monomials which are divisible by $Q_{3,1}^{\frac{q^{\ell+1}-1}{q-1}}$ or by $Q_{3,0}$. 
\end{proof}


 	
 The following corollary will be useful for induction. 
\begin{corollary}\label{cor: reduction}  Suppose $0 \leq \ell \leq m-2$ and $m-3 \geq \lambda_2 \geq \lambda_3 \geq 0$. The following equalities hold in $S$:  
\begin{enumerate}
    \item   
    \begin{multline*}
 \delta_{3;m-\ell-1}(1)^{q^{\ell}}  \delta_{3;m} (Q_{2,1}^{q^{\ell}} Q_{2,0}^{\frac{q^{\ell}-1}{q-1}}) - \delta_{3;m-\ell-1}(Q_{2,1})^{q^{\ell}} \delta_{3;m} ( Q_{2,0}^{\frac{q^{\ell}-1}{q-1}}) = Q_{3,1}^{\frac{q^{m-2}-q^{\ell}}{q-1}} Q_{3,0}^{\frac{q^{\ell}-1}{q-1}} + \text{other terms},
     \end{multline*}
   where the other terms are in the ideal $(Q_{3,0}^{\frac{q^{\ell+1}-1}{q-1}})$. 
   \item More generally, 
    \begin{multline*}
\delta_{3;\lambda_2-\lambda_3+1}(1)^{q^{\lambda_3}} \delta_{3;m} (Q_{2,1}^{q^{\lambda_3}}  Q_{2,0}^{\frac{q^{\lambda_3}-1}{q-1}}) - \delta_{3;\lambda_2-\lambda_3+1}(Q_{2,1})^{q^{\lambda_3}} \delta_{3;m} ( Q_{2,0}^{\frac{q^{\lambda_3}-1}{q-1}})  = \\  Q_{3,2}^{\frac{q^{m-2}-q^{\lambda_2}}{q-1}} Q_{3,1}^{\frac{q^{\lambda_2}-q^{\lambda_3}}{q-1}}  Q_{3,0}^{\frac{q^{\lambda_3}-1}{q-1}}  + \text{other terms},  
    \end{multline*}
    where the other terms are in the ideal $(Q_{3,0}^{\frac{q^{\lambda_3+1}-1}{q-1}}, Q_{3,1}^{\frac{q^{\lambda_2+1}-q^{\lambda_3}}{q-1}} Q_{3,0}^{\frac{q^{\lambda_3-1}}{q-1}})$.  
\end{enumerate}
\end{corollary}	 
\begin{proof}
This is just combining Lemma \ref{lem: delta3 formula} and Proposition \ref{prop:key reduction 3}. 

For the first part, note that 
\[
\delta_{3;m} (Q_{2,1}^{q^{\ell}} Q_{2,0}^{\frac{q^{\ell}-1}{q-1}}) = \delta_{3;m} (Q_{2,1}^{\frac{q^{\ell+1}-q^{\ell}}{q-1}} Q_{2,0}^{\frac{q^{\ell}-1}{q-1}})  = \delta_{3;m-\ell} (Q_{2,1})^{q^{\ell}} Q_{3,0}^{\frac{q^{\ell}-1}{q-1}},
\]
and similarly
\[
\delta_{3;m} ( Q_{2,0}^{\frac{q^{\ell}-1}{q-1}}) = \delta_{3;m-\ell}(1)^{q^{\ell}}Q_{3,0}^{\frac{q^{\ell}-1}{q-1}}. 
\]
The left-hand side of the required formula becomes
\[
\bigg(\delta_{3;m-\ell-1}(1) \delta_{3;m-\ell} (Q_{2,1}) - \delta_{3;m-\ell-1}(Q_{2,1}) \delta_{3;m-\ell} (1)\bigg)^{q^{\ell}} Q_{3,0}^{\frac{q^{\ell}-1}{q-1}}. 
\]
We can then apply \ref{prop:key reduction 3} part (i) to finish the proof.  

The second part is proved similarly.
\end{proof}

\begin{corollary}\label{cor:edge elements rank 3}
The Dickson monomials at the edge of $\Delta^m_3$ listed below are also in the span of $\mathcal{B}_{m}(3)$: 
		\begin{enumerate}
			\item $Q_{3,0}^{\frac{q^{m-2}-1}{q-1}}$,  
			\item $Q_{3,1}^{\frac{q^{m-2}-q^{\lambda_3}}{q-1}} Q_{3,0}^{\frac{q^{\lambda_3}-1}{q-1}} $, where $m-3 \geq \lambda_3$,
		  	\item $Q_{3,2}^{\frac{q^{m-2}-q^{\lambda_2}}{q-1}} Q_{3,1}^{\frac{q^{\lambda_2}-q^{\lambda_3}}{q-1}} Q_{3,0}^{\frac{q^{\lambda_3}-1}{q-1}}$, where $m-3 \geq \lambda_2 \geq \lambda_3$. 
		\end{enumerate}	    
\end{corollary}
\begin{proof}
We have already seen from part (1) of Lemma \ref{lem: delta3 formula} that the monomial $Q_{3,0}^{\frac{q^{m-2}-1}{q-1}}$ belongs to $\delta_3 (\Delta^m_2)$. For the remaining $2$ families, we will show that they belong to the span of $\mathcal{B}_m(3)$ by downward induction on $\lambda_3$.  
 
 When $\lambda_3 = m-3$, it follows that $\lambda_2 = \lambda_3 = m-3$. The edge monomials 
 \[
 Q_{3,1}^{q^{m-3}} Q_{3,0}^{\frac{q^{m-3}-1}{q-1}} \quad \text{and} \quad Q_{3,2}^{q^{m-3}} Q_{3,0}^{\frac{q^{m-3}-1}{q-1}}
 \]
 belong to the span of $\mathcal{B}_{m}(3)$ by parts (1) and (3) of Lemma \ref{lem: delta3 formula}. In addition, by Proposition \ref{prop: delta and Dickson}, any Dickson monomial divisible by these monomials is also in the image of $\delta_3$, and thus belongs to the span of $\mathcal{B}_m(3)$. 
 
 Now suppose, as induction hypothesis, that all edge monomials with $Q_{3,0}$-exponent $\frac{q^{\lambda_3+1}-1}{q-1}$ as well as their multiples are already in the span of $\mathcal{B}_m(3)$. We consider the edge monomials corresponding to $\lambda_3$:
 \begin{enumerate}
 	\item $Q_{3,1}^{\frac{q^{m-2}-q^{\lambda_3}}{q-1}} Q_{3,0}^{\frac{q^{\lambda_3}-1}{q-1}}$, 
   \item $Q_{3,2}^{\frac{q^{m-2}-q^{\lambda_2}}{q-1}} Q_{3,1}^{\frac{q^{\lambda_2}-q^{\lambda_3}}{q-1}} Q_{3,0}^{\frac{q^{\lambda_3}-1}{q-1}}$. 
 \end{enumerate}
 Part (i) of Corollary \ref{cor: reduction} shows that $Q_{3,1}^{\frac{q^{m-2}-q^{\lambda_3}}{q-1}} Q_{3,0}^{\frac{q^{\lambda_3}-1}{q-1}}$ lies in the image of $\delta_3 (\Delta^m_2)$ modulo terms with strictly higher $\lambda_3$. By the induction hypothesis, these higher terms are already in the span of $\mathcal{B}_{m}(3)$.

For the second monomial, observe that if there exists a monomial with the same  $Q_{3,0}$-exponent $\frac{q^{\lambda_3}-1}{q-1}$ but a higher exponent of $Q_{3,1}$, such that it becomes non-essential, then this monomial must be a multiple of 
\[
Q_{3,1}^{\frac{q^{m-2}-q^{\lambda_3}}{q-1}} Q_{3,0}^{\frac{q^{\lambda_3}-1}{q-1}},
\]
which by the previous step, is already in the span of $\mathcal{B}_{m}(3)$. Alternatively, if the $Q_{3,0}$-exponent is strictly greater than $\frac{q^{\lambda_3}-1}{q-1}$, the induction hypothesis can be applied directly. 
\end{proof}	

 Now we are ready to finish the proof of the claim that $\mathcal B_m (3) $ is a generating set for $\Q_m (3)^G $. Corollary \ref{cor:edge elements rank 3} shows that the subspace of $\Q_m(3)^G$ spanned by $\mathcal{B}_m(3)$ is closed under the action of the Dickson algebra $\D_3$. Indeed, we have already computed the action of $\D_3$ on $\delta_{s+1}^{3-s} (\Delta^m_s)$ for $0 \leq s \leq 2$ in  \ref{prop: delta and Dickson}, and the product of $Q_{3,i}$ with a monomial in $\Delta^m_3$ is clearly an edge monomial, which is in the span of $\mathcal{B}_m(3)$ by Corollary \ref{cor:edge elements rank 3}.   
 The proof is finished. 
 

     \section{A filtration of submodules over the Steenrod algebra}\label{sec: Steenrod}
 In this final section, we investigate the filtration $\mathcal{F}_{n,k}$ of $\Q_m(n)^G$ as defined in the Introduction \ref{def:filtration}.  Recall that for each positive integer $n$, $0 \leq k \leq \min (m,n)$,   $\mathcal{F}_{n,k}$  is the $\mathbb{F}_q$-subspace of $\Q_m(n)^G$:  
     \[
\mathcal{F}_{n,k} = \Span \{\delta_{s+1}^{n-s} (f) \colon f \in \Delta^m_s, 0 \leq s \leq \min(m,k) \}. 
     \]
The aim of this final section is to demonstrate Theorem \ref{thm:filtration}. For $n \leq 3$, we have shown that this increasing filtration of subspaces is exhaustive. We begin with the following  observation which allows more flexibility in working with $\mathcal{F}_{n,k}$: 
\begin{lemma}\label{def: filtration more flexible} Let $n \leq 3$. For each $0 \leq k \leq \min (m,n)$, we have 
 \[
\mathcal{F}_{n,k} = \Span \{\delta_{s+1}^{n-s} (f) \colon f \in \D_s, 0 \leq s \leq \min(m,k) \}. 
     \]
\end{lemma}
The point is that this filtration does not depend on the subspace $\Delta^m_s$.   
\begin{proof}
For $n=2$ this is essentially pointed out in Example \ref{y family} and the proof of Proposition \ref{prop:GL2}.  For $n=3$, this is a combination of Example \ref{a family} and the proof of Proposition \ref{prop:GL3}. 
\end{proof}
\begin{corollary}\label{cor: filtration of Dickson module}
    For $n \leq 3$, $\mathcal{F}_{n,*}$ is a filtration of $\mathcal{D}_n$-submodules of the invariant ring $\Q_m(n)^G$. Moreover, for $0 \leq k < \min (m,n)$, $\mathcal{F}_{n,k}$ is annihilated by $Q_{n,0}, \ldots, Q_{n,n-k-1}$. 
\end{corollary}
\begin{proof}
    This is immediate from the above Lemma and Proposition \ref{prop: delta and Dickson}. 
\end{proof}
 The calculations in the previous sections demonstrate that our basis $\mathcal{B}_m(3)$ for $\Q_m(3)^G$ is natural in the sense that the resulting filtration is a filtration by $\D_3$-submodules. Moreover, the action of the Dickson algebra follows a pattern similar to the one discussed in \cite{LewisReinerStanton2017}, where the structure of the cofixed space $S_G$ is described as an $S^G$-module. A similar pattern also emerges in the case $n=2$. 

 We now turn our attention to the Steenrod algebra. For background on the mod $q$ Steenrod algebra $\mathcal{A}$ and its role in modular invariant theory, we refer to Larry Smith's article \cite{LarrySmith2007}, which provides a purely algebraic treatment of the Steenrod algebra over a general finite field.    

Recall that for each $i \geq 0$, the $i$th Steenrod reduced power operation $\mathcal{P}^i$ is an $\mathbb{F}_q$-linear map satisfying the following conditions:
\begin{enumerate}
    \item The unstable condition: $\mathcal{P}^i (f) = f^q$ if $i = \deg (f)$ and $\mathcal{P}^i (f) = 0$ if $i > \deg (f)$. 
    \item The Cartan formulae: $\mathcal{P}^i (fg) = \sum_{a+b=i} \mathcal{P}^a (f) \mathcal{P}^b (g)$. 
\end{enumerate}
Moreover, $\mathcal{P}^0$ acts as the identity. Over the polynomial ring $S = \mathbb{F}_q [x_1, \ldots, x_n]$, in which each $x_i$ is given degree $1$, the action of the Steenrod algebra is completely determined by the unstable condition and the Cartan formulae. More concretely, we have   
\[
\mathcal{P}^i (v^j) = \binom{j}{i} v^{j + i(q-1)},
\]
for any linear polynomial $v$. Here the binomial coefficient is taken modulo $p$. Thus, if $j$ is a $q$-power, then $\mathcal{P}^i (v^j) \neq 0$ iff $i=0$ (in which case $\mathcal{P}^0$ is the identity operator) or $i=j$ (in which case $\mathcal{P}^i (v^i) = v^{qi}$ is the Frobenius operator). It is convenient to make the convention that $\mathcal{P}^k = 0$ if $k < 0$. 

The Steenrod algebra action and the group $G$ action on $S$ commute. Moreover, it is evident that the Frobenius ideal $I_m$ is stable under the action of the Steenrod operations. Consequently, the quotient $\Q_m = S/I_m$ inherits the structure of an $\mathcal A$-module. It follows that each Steenrod reduced power $\mathcal{P}^i$ induces a well-defined endomorphism of $\Q_m^{G}$. In other words, the subspace of invariants $\Q_m^{G}$ is closed under the action of the Steenrod algebra, and therefore forms an $\mathcal{A}$-module.

The proof of Theorem \ref{thm:filtration} uses several lemmas. We begin with a general formula describing how the Steenrod operations interact with the delta operator.  

\begin{lemma}\label{lem: Steenrod action on delta2} 
 If $f$ is a polynomial which is $G_1$-invariant in the first variable, then for any $k \geq 0$, the following holds in $S$:
\begin{multline*}
    \mathcal{P}^k (\delta_2 f)  + Q_{2,1} \mathcal{P}^{k-q} (\delta_2 f) + Q_{2,0} \mathcal{P}^{k-q-1} (\delta_2 f)   =\\ \delta_{2} (Q_{1,0} \mathcal{P}^{k-1} f) +  \delta_{2} (\mathcal{P}^{k} f) +\delta_{2;m+1} (Q_{1,0} \mathcal{P}^{k-1-q^m} f) + \delta_{2;m+1} (\mathcal{P}^{k-q^m} f). 
\end{multline*} 
\end{lemma}	
\begin{proof} The condition imposed on $f$ assures that $\delta_2(f)$ as well as the other terms in this lemma are polynomial.
The equation above is then obtained by applying the Cartan formula for the product $L_2 \delta_{2} (f)$ and then dividing both sides by $L_2$. Note that $\mathcal{P}^q L_2 = Q_{2,1} L_2$, $\mathcal{P}^{q+1} L_2 = Q_{2,0} L_2$ and $\mathcal{P}^k L_2 = 0$ for all other $k \geq 1$. 
\end{proof}
\begin{corollary}\label{cor:F21}
 $\mathcal{F}_{2,1}$ is a submodule of $\Q_m(2)^G$ over the Steenrod algebra. 
\end{corollary}
\begin{proof}
We must show that for each $0 \leq s <[m]_q$, $\mathcal{P}^k (\delta_2Q_{1,0}^s) \in \mathcal{F}_{2,1}$ for all $k > 0$. If $k \geq q^m$, the degree of $\mathcal{P}^k (\delta_2Q_{1,0}^s)$ is at least 
\[
q^m (q-1) + q^m-q + s(q-1) \geq 2(q^m-1). 
\]
So we can assume that $k < q^m$. In this case, when we apply the formula of Lemma \ref{lem: Steenrod action on delta2} for $f = Q_{1,0}^a$,  the last two summands on the right-hand side vanish. Thus we have an equality in $S$: 
\begin{equation}
     \mathcal{P}^k (\delta_2 Q_{1,0}^s)  + Q_{2,1} \mathcal{P}^{k-q} (\delta_2 Q_{1,0}^s) + Q_{2,0} \mathcal{P}^{k-q-1} (\delta_2 Q_{1,0}^s)   = \delta_{2} (Q_{1,0} \mathcal{P}^{k-1} Q_{1,0}^s) +  \delta_{2} (\mathcal{P}^{k} Q_{1,0}^s). 
\end{equation}
Both summands on the right-hand side are evidently of the form $\delta_2 (Q_{1,0}^s)$. Now we can proceed by induction on $k$ using the above formula and the calculation in Proposition \ref{prop: delta and Dickson} when working in $\Q_m(2)$.    
\end{proof}
 Note that $\mathcal{F}_{2,0}$ is spanned by the top class $x_1^{q^m-1}x_2^{q^m-1}$, and $\mathcal{F}_{2,2} = \Q_m(2)^G$. Both are evidently $\mathcal{A}$-submodules of $\Q_m(2)^G$.

\begin{proposition}\label{prop: F31}
$\mathcal{F}_{3,1}$ is a submodule of $\Q_m(3)^G$ over the Steenrod algebra. 
\end{proposition}
\begin{proof}
We need to show that for each $k\geq 1$ and $0 \leq a < [m]_q$, $\mathcal{P}^k (\delta_2^2 Q_{1,0}^s)$ is in the span of $\delta_1^3 (1)$ and $\delta_2^2 (\Delta^m_1)$. Note that $\delta_1^3 (1)$ spans the one-dimensional space of $\Q_m(3)^G$ in the top degree $3(q^m-1)$. 

First of all, if $k \geq q^m$, then the degree of $\mathcal{P}^k (\delta_2^2 Q_{1,0}^s)$, which is 
$q^m (q-1) + 2(q^m-q) + s(q-1)$ is at least $3(q^m-1)$ unless $q=2$ and $k=2^m$, $s = 0$. In this range, 
$\Q_m(3)$ is either spanned by $\delta_3^3 (1)$ or zero. 

If $q=2$, $k=2^m$, $s=0$, then $\mathcal{P}^{2^m} (\delta_2^2 (1))$ is in degree $3(2^m-1)-1$. We claim that there is no $G$-invariant polynomial in this degree. Indeed, the only non-zero symmetric polynomial in this degree is 
\[
x_1^{2^m-1} x_2^{2^m-1} x_3^{2^m-2} + x_1^{2^m-1} x_2^{2^m-2} x_3^{2^m-1} + x_1^{2^m-2} x_2^{2^m-1} x_3^{2^m-2}.
\]
By direct inspection, this polynomial is not invariant under the action of the operation mapping $x_1 \mapsto x_1 + x_2$ and fixes $x_2, x_3$.  

Now let us consider the case $k < q^m$. By Lemma \ref{lem: Steenrod action on delta2}, we have the following equality in $S$: 
\begin{equation}\label{eq: 8.5 simplified}
     \mathcal{P}^k (\delta_2 f)  + Q_{2,1} \mathcal{P}^{k-q} (\delta_2 f) + Q_{2,0} \mathcal{P}^{k-q-1} (\delta_2 f)   = \delta_{2} (Q_{1,0} \mathcal{P}^{k-1} f) +  \delta_{2} (\mathcal{P}^{k} f)
\end{equation}
 for all polynomial $f$ which is $GL_1$-invariant in the first variable. Applying this for $f = Q_{1,0}^s$, it follows by induction on $k$ that for $0 \leq k < q^m$, $\mathcal{P}^k (\delta_2 Q_{1,0}^s)$ is a linear combination of polynomials of the form 
\[
Q_{2,1}^{a} Q_{2,0}^{b} \delta_2 (Q_{1,0}^c).
\]
Now applying \eqref{eq: 8.5 simplified} again but for $f = \delta_2 (Q_{1,0}^s)$, we get 
\begin{equation}\label{eq: 8.5 simplified bis}
     \mathcal{P}^k (\delta_2^2(Q_{1,0}^s))  + Q_{2,1} \mathcal{P}^{k-q} (\delta_2^2 (Q_{1,0}^s)) + Q_{2,0} \mathcal{P}^{k-q-1} (\delta_2^2 (Q_{1,0}^s))   = \delta_{2} (Q_{1,0} \mathcal{P}^{k-1} \delta_2(Q_{1,0}^s)) +  \delta_{2} (\mathcal{P}^{k} \delta_2(Q_{1,0}^s)).
\end{equation}
The right-hand side is then a sum of polynomials of the following forms
\begin{itemize}
    \item $\delta_2 (Q_{1,0} Q_{2,1}^{a} Q_{2,0}^{b} \delta_2 (Q_{1,0}^c))$, or 
    \item $\delta_2 (Q_{2,1}^{a} Q_{2,0}^{b} \delta_2 (Q_{1,0}^c))$.
\end{itemize}
  Up to this point, we have been working in $S$. Now we consider the situation in $\Q_m(3)$. From Lemma \ref{lem: delta2 f delta2}, we see that this sum belongs to the subspace of $\Q_m(3)$ spanned by polynomials of the forms
  \begin{itemize}
      \item $\delta_2^2 (Q_{1,0}^a)$, or
      \item $x_1^{q^m-1} x_2^{q^m-1} x_3^{a(q-1)}$. 
  \end{itemize}
Both sides of Equation \ref{eq: 8.5 simplified bis} now are $P(2,1)$-invariant. Take the transfer from $P(2,1)$ to $G$. We have
$\mathcal{P}^i (\delta_2^2 (Q_{1,0}^s))$ is $G$-invariant in $\Q_m(3)$ and 
\[
\tr_{P(2,1)}^G (Q_{2,1}) = \tr_{P(2,1)}^G (Q_{2,0}) = \tr_{P(2,1)}^{G} (x_1^{q^m-1} x_2^{q^m-1} x_3^{a(q-1))}) = 0.  
\]
The first two terms vanish in $S$ for degree reasons: there are no genuine $G$-invariants in degrees $q^2-q$ and $q^2-1$. The third vanishes in $\Q_m(3)$ because if $x_1^{q^m-1} x_2^{q^m-1} x_3^{a(q-1)}$ is a nontrivial summand of a $G$-invariant polynomial in $\Q_m(3)$, then it is a $P(1,2)$-invariant as well, but from Remark \ref{rem:G12}, this can only happen at the top degree. 

 It follows that $\mathcal{P}^k (\delta_2^2 (Q_{1,0}^s))$ must be a scalar multiple of some $\delta_2^2 (Q_{1,0}^a)$, and thus belongs to $\mathcal{F}_{3,1}$. 
\end{proof}

We now consider the interaction between the Steenrod powers and the operator $\delta_3$.

\begin{lemma}\label{lem:Steenrod action on delta3} 
If $f$ is a polynomial which is $G_2$-invariant in the first two variables, then for any $k \geq 0$, the following holds in $S$:
 \begin{multline*}
	 	\mathcal{P}^k (\delta_3 f)  + Q_{3,2} \mathcal{P}^{k-q^2} (\delta_3 f) + Q_{3,1} \mathcal{P}^{k-q^2-q} (\delta_3 f) +Q_{3,0} \mathcal{P}^{k-q^2-q-1} (\delta_3 f)    =\\ 
	 	\delta_3\bigg( \mathcal P^k (f)+ Q_{2,1} \mathcal P^{k-q}(f) +Q_{2,0}\mathcal P^{k-q-1} (f)\bigg) + 	\delta_{3;m+1}\bigg( \mathcal P^{k-q^m} (f)+ Q_{2,1} \mathcal P^{k-q^m-q}(f) +Q_{2,0}\mathcal P^{k-q^m-q-1} (f)\bigg).
	 \end{multline*} 
\end{lemma}			
\begin{proof}
The above equality is again an application of the Cartan formula to $\mathcal{P}^k (L_3 \delta_3 f)$, using the formulae for the action of the Steenrod powers on $L_2$ given in the proof of Lemma \ref{lem: Steenrod action on delta2}, and similar formulae for these actions on $L_3$:
$\mathcal{P}^{q^2} L_3 = Q_{3,2} L_3$, $\mathcal{P}^{q^2+q} L_3 = Q_{3,1} L_3$, $\mathcal{P}^{q^2+q+1} L_3 = Q_{3,0} L_3$ and $\mathcal{P}^k L_3 = 0$ for all other $k \geq 1$. 
\end{proof}

The next lemma discusses how to deal with the last term in the above formula.

\begin{lemma}\label{lem:delta_{3,m+1}} If $h \in \D_2$, then the following holds in $\Q_m$, 
\[
\delta_{3;m+1} (h) = 
\begin{cases}
 0 & \text{if $q> 3$ or $q=3$ and $\deg (h) > 0$;}\\
 \delta_{3;m} (Q_{2,1}^{q^{m-1}}) & \text{if $q= 3$ and $h=1$;}\\
 \in \Span \{\delta_{2;m}^2 (\Delta^m_1), \delta_{1;m}^3 (1)\} & \text{if $q=2$.}\\
\end{cases}
\]
\end{lemma}
\begin{proof}
Recall from the fundamental equation that $X^{q^3} - Q_{3,2} X^{q^2} + Q_{3,1} X^q - Q_{3,0} X = 0$ if $X$ is in the span of $x_1, x_2, x_3$. It follows that 
\[
x_i^{q^{m+1}} - Q_{3,2}^{q^{m-2}} x_i^{q^m} + Q_{3,1}^{q^{m-2}} x_i^{q^{m-1}} - Q_{3,0}^{q^{m-2}} x_i^{q^{m-2}} = 0, \quad 1 \leq i \leq 3,
\]
and we obtain the following identity in $S$ relating $\delta_{3}$ for various $m$: 
\[
\delta_{3;m+1} (h) - Q_{3,2}^{q^{m-2}} \delta_{3;m} (h) + Q_{3,1}^{q^{m-2}} \delta_{3;m-1} (h) - Q_{3,0}^{q^{m-2}} \delta_{3;m-2} (h) = 0.
\]
Note that the degree of $\delta_{3;m+1} (h)$ equals $q^{m+1}-q^2 + \deg (h)$. If $q>3$, this is strictly greater than $3(q^m-1)$, and we can conclude immediately that $\delta_{3;m+1} (h) = 0$ in $\mathcal{Q}_m(3)$.  

Consider now the case $q=3$. Observe that, by direct computation, $Q_{3,0} = Q_{3,1} = 0$, $Q_{3,2} = (x_1 x_2 x_3)^{3^2-3^1}$ in $\Q_2(3)$. Hence $Q_{3,1}^{q^{m-2}} = Q_{3,0}^{q^{m-2}} = 0$ in $\Q_m(3)$ and the above equation becomes
\begin{equation}\label{eq:delta3m+1}
    \delta_{3;m+1}(h) = Q_{3,2}^{q^{m-2}} \delta_{3;m} (h) = (x_1 x_2 x_3)^{3^m-3^{m-1}} \delta_{3;m} (h) \quad \text{in $\Q_m(3)$}. 
\end{equation}
 On the other hand, for degree reasons, we must have $\deg (h) \leq 3(3^m-1)-(3^{m+1}-3^2) = 6$. In this range, there are (up to a scalar) only two Dickson polynomials, namely $h=1$ or $h=Q_{2,1}$.
If $h=1$, then from Proposition \ref{prop: delta and Dickson}, we have immediately  
\[
\delta_{3;m+1}(1) = Q_{3,2}^{3^{m-2}} \delta_{3;m} (1) = \delta_{3;m}(Q_{2,1}^{3^{m-1}}). 
\]
If $h=Q_{2,1}$, we prove by induction on $m \geq 2$ that 
\[
\delta_{3;m+1} (Q_{2,1}) = 0 \quad \text{in $\Q_m(3)$}. 
\]
Indeed, the base case can be checked directly since $\delta_{3;3} (Q_{2,1}) = Q_{3,1}$ in $S$, and the induction step follows immediately from Equation \eqref{eq:delta3m+1}.  


The next and final case is $q=2$. Using the identity 
\[
x_1^{2^{m+1}} = x_1^{2^m} (x_1^{2^m} + x_2^{2^m} + x_3^{2^m}) - x_1^{2^m} (x_2^{2^m} + x_3^{2^m}), 
\]
and similarly for $x_2^{2^{m+1}}$ and $x_3^{2^{m+1}}$, we see that $\delta_{3;m+1} (h)$ equals in $\Q_m(3)$ to 
\[
\frac{1}{L_3} \begin{vmatrix}
    x_1 & x_2 & x_3\\
    x_1^2 & x_2^2 & x_3^2\\
    x_1^{2^m} h(x_2,x_3)(x_2^{2^m} + x_3^{2^m}) & x_2^{2^m} h(x_1,x_3)(x_1^{2^m} + x_3^{2^m}) & x_3^{2^m} h(x_1,x_2)(x_1^{2^m} + x_2^{2^m})  
\end{vmatrix}.
\]
It is straightforward to verify the following equality in $S$:  
\[
x_1^{2^m} + x_2^{2^m} = Q_{2,1} \delta_{2;m} (1) - \delta_{2;m} (Q_{1,0}^2),
\]
We then obtain an equality in $\Q_m(3)$: 
\[
\delta_{3;m+1}(h) = \delta_{3;m} (Q_{2,1}h \delta_{2;m}(1)) - \delta_{3;m} (h \delta_{2;m}(Q_{1,0}^2)). 
\]
The conclusion follows from Corollary \ref{cor 73}.   
\end{proof}

\begin{remark}
With a little more work, using \ref{lem: delta2 f delta2}, it can be shown that when $q=2$, then for $h \in \D_2$, we have the following explicit formula for $\delta_{3;m+1} (h)$ in $\Q_m(3)$:  
\[
\delta_{3;m+1} (h) = 
\begin{cases}
0  & \text{if $h$ is a multiple of $Q_{2,0}$;}\\
\delta_{2;m}^2 (Q_{1,0}^{2s}) & \text{if $h = Q_{2,1}^s$, $s \geq 0$}.
\end{cases}
\]
We leave the details for the interested reader. 
\end{remark}

Now we can finish the proof of Theorem \ref{thm:filtration}. It remains to show that $\mathcal{F}_{3,2}$ is an $\mathcal{A}$-submodule of $\Q_m(3)^G$. From Lemma \ref{lem:delta_{3,m+1}}, we see that if $f \in \D_2$, then 
\[
\mathcal{P}^k (\delta_3 f)  + Q_{3,2} \mathcal{P}^{k-q^2} (\delta_3 f) + Q_{3,1} \mathcal{P}^{k-q^2-q} (\delta_3 f) +Q_{3,0} \mathcal{P}^{k-q^2-q-1} (\delta_3 f)   \in \mathcal F_{3,2}.
\]
We can proceed by induction. If $\mathcal{P}^i (\delta_3 f) \in \mathcal{F}_{3,2}$ for all $i<k$, then 
 $Q_{3,j} \mathcal{P}^i (\delta_3 f)$ also belongs to $\mathcal{F}_{3,2}$ by Proposition \ref{prop: delta and Dickson}. It follows that $\mathcal{P}^k (\delta_3 f)$ also belongs to $\mathcal F_{3,2}$ and the proof is complete.

			\bibliographystyle{amsplain}
			
			\providecommand{\bysame}{\leavevmode\hbox to3em{\hrulefill}\thinspace}
			\providecommand{\MR}{\relax\ifhmode\unskip\space\fi MR }
			\providecommand{\MRhref}[2]{%
				\href{http://www.ams.org/mathscinet-getitem?mr=#1}{#2}
			}
			\providecommand{\href}[2]{#2}

			

			
		\end{document}